\documentclass[final]{siamltex}
\usepackage{graphicx}  
\usepackage{psfrag}  
\usepackage{amsmath}  
\usepackage{amssymb}  
\usepackage{amsfonts}  

\renewcommand{\vec}[1]{\mbox{\boldmath $#1$}}  
\newcommand{\R}{\mathbb{R}}  
\newcommand{\C}{\mathbb{C}}  
\newcommand{\N}{\mathbb{N}}  
\newcommand{\curl}{\mbox{\rm curl}\; }  
\newtheorem{remark}{Remark}  

\begin{document} 

\title{Optimized Schwarz Methods for Maxwell's equations}

\author{V. Dolean\footnote{Univ. de Nice Sophia-Antipolis, Laboratoire
J.-A. Dieudonn\'e, Nice, France.  dolean@unice.fr},
M.J. Gander\footnote{Section de Math\'ematiques, Universit\'e de
Gen\`eve, CP 64, 1211 Gen\`eve, Martin.Gander@math.unige.ch} \and
L. Gerardo-Giorda\footnote{Department. of Mathematics, University of
Trento, Italy. gerardo@science.unitn.it}}

\maketitle  
    
\begin{abstract}  
  Over the last two decades, classical Schwarz methods have been
  extended to systems of hyperbolic partial differential equations,
  using characteristic transmission conditions, and it has been
  observed that the classical Schwarz method can be convergent even
  without overlap in certain cases. This is in strong contrast to the
  behavior of classical Schwarz methods applied to elliptic problems,
  for which overlap is essential for convergence.  More recently,
  optimized Schwarz methods have been developed for elliptic partial
  differential equations. These methods use more effective
  transmission conditions between subdomains than the classical
  Dirichlet conditions, and optimized Schwarz methods can be used both
  with and without overlap for elliptic problems. A simple computation
  shows why the classical Schwarz method applied to both the time
  harmonic and time discretized Maxwell's equations converges without
  overlap: for a given frequency we obtain the same convergence rate
  as for an optimized Schwarz method for a scalar
  elliptic equation. Based on this insight, we show how to develop an
  entire new hierarchy of optimized overlapping and non-overlapping
  Schwarz methods for Maxwell's equations with greatly enhanced
  performance compared to the classical Schwarz method. We also derive
  for each algorithm asymptotic formulas for the optimized
  transmission conditions, which can easily be used in implementations
  of the algorithms for problems with variable coefficients. We
  illustrate our findings with numerical experiments.
\end{abstract} 

  

\pagestyle{myheadings} 
\thispagestyle{plain} 
\markboth{V. DOLEAN
M.J. GANDER AND L. GERARDO-GIORDA}{OPTIMIZED SCHWARZ METHODS FOR
MAXWELL'S EQUATIONS}

  
\section{Introduction}\label{intro}  
Schwarz algorithms have experienced a second youth over the last
decades, when distributed computers became more and more powerful and
available. Fundamental convergence results for the classical Schwarz
methods were derived for many partial differential equations, and can
now be found in several authoritative reviews, see
\cite{Chan:1994:DOA,Xu:1992:IMS,Xu:1998:NDD}, and books, see
\cite{Smith:1996:DPM,Quarteroni:1999:DDM,Toselli:2004:DDM}.  The
Schwarz methods were also extended to systems of partial differential
equations, such as the time harmonic Maxwell's equations, see
\cite{Despres:1992:ADD,Collino:1997:NIC}, or the time discretized
Maxwell's equations, see \cite{Toselli:2000:OSM}, or to linear
elasticity \cite{Faccioli:1996:SDM,Faccioli:1997:EWP}, but much less
is known about the behavior of the Schwarz methods applied to
hyperbolic systems of equations. This is true in particular for the
Euler equations, to which the Schwarz algorithm was first applied in
\cite{Quarteroni:1990:DCL,Quarteroni:1996:HDD}, where classical
(characteristic) transmission conditions are used at the interfaces,
or with more general transmission conditions in \cite{Clerc:1998:NOS}.
The analysis of such algorithms applied to systems proved to be very
different from the scalar case, see
\cite{Dolean:2002:CIC,Dolean:2004:CAS}.
  
Over the last decade, a new class of overlapping Schwarz methods was
developed for scalar partial differential equations, namely the
optimized Schwarz methods. These methods are based on a classical
overlapping domain decomposition, but they use more effective
transmission conditions than the classical Dirichlet conditions at the
interfaces between subdomains. New transmission conditions were
originally proposed for three different reasons: first, to obtain
Schwarz algorithms that are convergent without overlap, see
\cite{Lions:1990:SAM} for Robin conditions. The second motivation for
changing the transmission conditions was to obtain a convergent
Schwarz method for the Helmholtz equation, where the classical
overlapping Schwarz algorithm is not convergent. As a remedy,
approximate radiation conditions were introduced in
\cite{Despres:1990:DDP,Despres:1992:ADD}.  The third motivation was
that the convergence rate of the classical Schwarz method is rather
slow and too strongly dependent on the size of the overlap. In a short
note on non-linear problems \cite{Hagstrom:1988:NED}, Hagstrom et
al. introduced Robin transmission conditions between subdomains and
suggested nonlocal operators for best performance. In
\cite{Charton:1991:MDD}, these optimal, non-local transmission
conditions were developed for advection-diffusion problems, with local
approximations for small viscosity, and low order frequency
approximations were proposed in \cite{Nataf:1995:FCD,Deng:1997:AND}.
In \cite{Sofronov:1998:NIO}, one can find low-frequency approximations
of absorbing boundary conditions for Euler equations. Independently,
at the algebraic level, generalized coupling conditions were
introduced in \cite{Tang:1992:GSS, Sun:1996:OSA} for discrete
overlapping Schwarz methods. Optimized transmission conditions for the
best performance of the Schwarz algorithm in a given class of local
transmission conditions were first introduced for advection diffusion
problems in \cite{Japhet:2000:OO2}, for the Helmholtz equation in
\cite{Chevalier:1998:SMO,Gander:2001:OSH}, and for Laplace's equation
in \cite{Enquist:1998:ABC}. For complete results and attainable
performance for a symmetric, positive definite problem, see
\cite{Gander:2006:OSM}, and for time dependant problems, see
\cite{Gander:2003:OSWW,Gander:2003:MRO}. The purpose of this paper is
to design and analyze a family of optimized overlapping and
non-overlapping Schwarz methods for Maxwell's equations, both for the
case of time discretized and time harmonic problems, and to provide
explicit formulas for the optimized parameters in the transmission
conditions of each algorithm in the family. These formulas can then
easily be used in implementations for Maxwell's equations with
variable coefficients. As we will see, one member of this family
reduces in the case of no overlap and constant coefficients to an
algorithm in a curl-curl formulation of Maxwell's equations, proposed
in \cite{Alonso:2006:NSM} based on \cite{Chevalier:1998:MNT}, which
already greatly enhanced the performance compared to the classical
approaches in \cite{Despres:1992:ADD,Collino:1997:NIC}.

This paper is organized as follows: in Section \ref{MaxwellSystem}, we
present Maxwell's equations and a reformulation thereof with
characteristic variables used in our analysis. In Section
\ref{TimeHarmonicSection}, we treat the case of time harmonic
solutions. We show that the classical Schwarz method for Maxwell's equations,
which uses characteristic Dirichlet transmission conditions between
subdomains is convergent even without overlap. Exploiting a parallel
with an optimized Schwarz method applied to an Helmholtz equation allows
us to develop an entirely new hierarchy of optimized Schwarz methods for
Maxwell's equations with greatly enhanced performance, both with and
without overlap. Similar equivalence has been presented in \cite{Dolean:2007:WCS} for the
Cauchy-Riemann equations.
In Section \ref{TimeDiscretizedSection}, we present and
analyze the corresponding hierarchy of optimized Schwarz methods for
time discretizations of Maxwell's equations. We then
show in Section \ref{NumericalSection} numerical experiments in two
and three spatial dimensions, both for the time harmonic and time
discretized case, which illustrate the performance of the new
optimized Schwarz methods for Maxwell's equations. We also include as
an application the cooking of a chicken in a microwave oven, a problem
with variable coefficients. In Section \ref{Conclusions}, we summarize
our findings and conclude with an outlook on future research
directions.  

\section{Maxwell's Equations}\label{MaxwellSystem}  
  
The hyperbolic system of Maxwell's equations describes the propagation
of electromagnetic waves. It is given by
\begin{equation} \label{maxw}   
  -\varepsilon \frac{\partial \vec{\cal E}}{\partial t} + \curl \vec{\cal H} -\sigma \vec{\cal E}= \vec{J},   
   \qquad   
  \mu \frac{\partial \vec{\cal H}}{\partial t} + \curl \vec{\cal E} = 0,    
\end{equation}   
where $\vec{\cal E}=({\cal E}_1,{\cal E}_2,{\cal E}_3)^T$ and
$\vec{\cal H}=({\cal H}_1,{\cal H}_2,{\cal H}_3)^T$ denote the
electric and magnetic fields, respectively, $\varepsilon$ is the {\it
electric permittivity}, $\mu$ is the {\it magnetic permeability},
$\sigma$ is the {\it electric conductivity} and $\vec{J}$ is the
applied current density. We assume the applied current density to be
divergence free, that is $\mbox{div} \vec{J}=0$. Denoting the vector of physical unknowns by
\begin{equation}\label{U=(E,H)}  
  \vec{u} = \left({\cal E}_1 , {\cal E}_2 , {\cal E}_3 , {\cal H}_1 , {\cal H}_2 , {\cal H}_3 \right)^T,  
\end{equation}  
Maxwell's equations (\ref{maxw}) can be rewritten in the form
\begin{equation}\label{maxwsyst}  
  (G+G_0\partial_t) \vec{u}+  G_x\partial_{x}\vec{u} + G_y \partial_{y} \vec{u}  
   + G_z\partial_{z} \vec{u} = (\vec{J};\vec{0}),  
\end{equation}  
where the coefficient matrices are  
$$  
G = \left[\begin{array}{cc}  
     \sigma I_3     &  \\  
     &  0_3  
   \end{array}\right],\,G_0 = \left[\begin{array}{cc}  
     \varepsilon I_3     &  \\  
     &  \mu I_3 
   \end{array}\right],\,  G_l  = \left[\begin{array}{cc}  
          & N_l \\  
     -N_l &    
   \end{array}\right],\qquad l=x,y,z,  
$$ 
where $0_3$ (resp. $I_3$) represent the $3\times 3$ zero (identity)
matrix, and the matrices $N_l$, $l=x,y,z$ are given by
$$  
  N_x=\left[\begin{array}{ccc}  
    0 & 0 & 0 \\  
    0 & 0 & 1 \\  
    0 & -1& 0  
  \end{array}\right],\qquad  
  N_y=\left[\begin{array}{ccc}  
    0 & 0 & -1 \\  
    0 & 0 & 0\\  
    1 & 0& 0  
  \end{array}\right],\qquad  
  N_z=\left[\begin{array}{ccc}  
    0 & 1 & 0 \\  
   -1 & 0 & 0\\  
    0 & 0 & 0  
  \end{array}\right].  
$$  
For any unit vector $\vec{n}=(n_1,n_2,n_3)$, $\|\vec{n}\|=1$, we can  
define the characteristic matrix of system (\ref{maxwsyst}) by  
$$ 
  C(\vec{n})  =  G_0^{-1}\left(n_1\left[\begin{array}{cc}  
       & N_x \\  
  -N_x &    
  \end{array}\right]+ n_2\left[\begin{array}{cc}  
       & N_y \\  
  -N_y &    
  \end{array}\right]+n_3 \left[\begin{array}{cc}  
       & N_z \\  
  -N_z &    
  \end{array}\right]\right),
$$   
whose eigenvalues are the characteristic speed of propagation along
the direction $\vec{n}$. A direct calculation shows that the matrix
$C(\vec{n})$ has real eigenvalues,
$$  
  \lambda_{1,2}=-c, \qquad \lambda_{3,4}=0, \qquad \lambda_{5,6}=c,  
$$ 
with $c = \frac{1}{\sqrt{\varepsilon\mu}}$ being the wave speed. This
implies that Maxwell's equations are hyperbolic, since the eigenvalues
are real, but not strictly hyperbolic, since the eigenvalues are not
distinct, see \cite{Serre:2007:MDH}. For the special case of the
normal vector $\vec{n}=(1,0,0)$, which we will use extensively later,
we obtain
$$  
  C(\vec{n})=\left(\begin{array}{cc}  
      & \frac{1}{\varepsilon}N_x \\  
  -\frac{1}{\mu}N_x &    
  \end{array}\right),
$$  
whose matrix of eigenvectors is given by
$$  
  L= \left[\begin{array}{cccccc}  
   0 & 0 & 0 & 1 & 0 & 0 \\  
  -{Z} & 0 & 0 & 0 & {Z} & 0 \\  
   0 & {Z} & 0 & 0 & 0 & -{Z} \\  
   0 & 0 & 1 & 0 & 0 & 0 \\  
   0 & 1 & 0 & 0 & 0 & 1 \\  
   1 & 0 & 0 & 0 & 1 & 0  
  \end{array}\right],  
$$ 
where $Z = \sqrt{\frac{\mu}{\varepsilon}}$ denotes the impedance. This
leads to the characteristic variables $\vec{w}=
(w_1,w_2,w_3,w_4,w_5,w_6)^T = L^{-1}\vec{u}$ associated with the
direction $\vec{n}$, where
\begin{equation}\label{charact}  
  \begin{array}{lll}  
   w_1= -\frac{1}{2}(\frac{1}{Z}{\cal E}_2-{\cal H}_3), \quad & w_2=
     \frac{1}{2}(\frac{1}{Z}{\cal E}_3+{\cal H}_2), \quad & w_3=
     {\cal H}_1, \\ 
     w_4= {\cal E}_1, 
     & w_5= \frac{1}{2}(\frac{1}{Z}{\cal E}_2+{\cal H}_3),& w_6=
     -\frac{1}{2}(\frac{1}{Z}{\cal E}_3-{\cal H}_2).
  \end{array}  
\end{equation}  
In the following, we will denote by $\vec{w}_+$, $\vec{w}_0$ and  
$\vec{w}_-$ the characteristic variables associated with the negative,  
zero, and positive eigenvalues respectively, that is  
\begin{equation}\label{W-W+}  
  \vec{w}_-=(w_1,w_2)^T, \qquad \vec{w}_0=(w_3,w_4)^T, \qquad  
    \vec{w}_+=(w_5,w_6)^T .  
\end{equation}  
Imposing classical or characteristic boundary conditions on a boundary
with unit outward normal vector $\vec{n}=(1,0,0)$ means to impose
Dirichlet conditions on the incoming characteristic variables
$\vec{w}_-$.  for a general normal vector $\vec{n}$, this is
equivalent to imposing the impedance condition (see
\cite{Serre:2007:MDH})
\begin{equation}\label{ImpedanceCondition}
  {\cal B}_{\vec{n}}(\vec{\cal E},\vec{\cal H}):= \vec{n}\times\frac{\vec{\cal
   E}}{Z}+\vec{n}\times(\vec{\cal H}\times\vec{n}) = \vec{s}.
\end{equation}

\section{Time Harmonic Solutions}\label{TimeHarmonicSection}   
  
Time harmonic solutions of Maxwell's equations are complex valued
static vector fields $\vec{E}$ and $\vec{H}$ such that the dynamic
fields
$$   
  \vec{\cal E}(\vec{x},t)={\cal R}e(\vec{E}(\vec{x})\exp(i \omega t)),\qquad  
  \vec{\cal H}(\vec{x},t)={\cal R}e(\vec{H}(\vec{x})\exp(i \omega t))  
$$   
satisfy Maxwell's equations (\ref{maxw}). The positive real parameter  
$\omega$ is called the {\em pulsation} of the harmonic wave. The harmonic  
solutions $\vec{E}$ and $\vec{H}$ satisfy the time-harmonic equations  
\begin{equation} \label{eq:THMaxwell} 
  - i \omega\varepsilon \vec{E}+ \curl \vec{H} -\sigma\vec{E} = \vec{J},
  \qquad i \omega\mu \vec{H}+\curl \vec{E}=\vec{0}.
\end{equation}  

\subsection{Classical and Optimized Schwarz Algorithm}
  \label{convTHTD}  
 
\begin{figure}
  \centering
  \psfrag{x}[][]{\footnotesize $x$}
  \psfrag{y}[][]{\footnotesize $y$}
  \psfrag{z}[][]{\footnotesize $z$}
  \psfrag{O1}[][]{\footnotesize $\Omega_1$}
  \psfrag{O2}[][]{\footnotesize $\Omega_2$}
  \psfrag{G12}[][]{\footnotesize $\Gamma_{12}$}
  \psfrag{G21}[][]{\footnotesize $\Gamma_{21}$}
  \psfrag{dO}[][]{\footnotesize $\partial \Omega$}
  \includegraphics[width=0.6\textwidth]{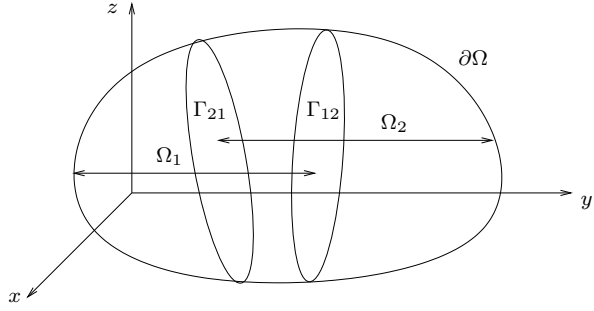}
  \caption{Overlapping domain decomposition.}
  \label{domaindec}
\end{figure}

We consider now the problem (\ref{eq:THMaxwell}) in a bounded domain
$\Omega$, with either Dirichlet conditions on the tangent electric
field, or impedance conditions, on $\partial \Omega$, in order to
obtain a well posed problem, see \cite{Nedelec:2001:AEE}.  In order to
explain the classical Schwarz algorithm for Maxwell's equation, we
decompose the domain into two overlapping subdomains $\Omega_1$ and
$\Omega_2$, as illustrated in Figure \ref{domaindec}. The
generalization of the algorithm formulation to the case of many
subdomains does not present any difficulties.
The classical Schwarz algorithm then solves for $n=1,2\ldots$ the
subdomain problems
\begin{equation}\label{ClassicalSchwarz}
  \begin{array}{rcll}  
    -i\omega\varepsilon\vec{E}^{1,n}
      +\curl\vec{H}^{1,n}-\sigma\vec{E}^{1,n} &=& \vec{J}
    \quad & \mbox{in $\Omega_1$} \\
  \qquad i \omega\mu \vec{H}^{1,n}+\curl \vec{E}^{1,n} & = & \vec{0}
    & \mbox{in $\Omega_1$} \\
  {\cal B}_{\vec{n}_1}(\vec{E}^{1,n},\vec{H}^{1,n})& = & 
  {\cal B}_{\vec{n}_1}(\vec{E}^{2,n-1},\vec{H}^{2,n-1})
    & \mbox{on $\Gamma_{12}$}\\
  - i \omega\varepsilon \vec{E}^{2,n}
      +\curl \vec{H}^{2,n} -\sigma\vec{E}^{2,n} & = & \vec{J}
    & \mbox{in $\Omega_2$} \\
  \qquad i \omega\mu \vec{H}+\curl \vec{E} & = & \vec{0}
    & \mbox{in $\Omega_2$} \\
  {\cal B}_{\vec{n}_2}(\vec{E}^{2,n},\vec{H}^{2,n})
    & = & 
  {\cal B}_{\vec{n}_2}(\vec{E}^{1,n-1},\vec{H}^{1,n-1})
    & \mbox{on $\Gamma_{21}$,} 
  \end{array}
\end{equation}
where $\Gamma_{12} = \partial\Omega_1\cap\Omega_2$, $\Gamma_{21} =
\partial\Omega_2\cap\Omega_1$ and ${\cal B}_{\vec{n}_j}$, $j=1,2$,
denotes the impedance boundary conditions defined in
(\ref{ImpedanceCondition}).  On the physical part of the boundary, the
given boundary conditions are imposed. While the choice of
transmission conditions ${\cal B}_{\vec{n}_j}$ is natural in the view
of the hyperbolic nature of the problem, we will see in our analysis
that there are better choices for the performance of the
algorithm. This leads to the so called optimized Schwarz methods,
\begin{equation}\label{OptimizedSchwarz}
  \begin{array}{rcll}  
    -i\omega\varepsilon\vec{E}^{1,n}
      +\curl\vec{H}^{1,n}-\sigma\vec{E}^{1,n} &=& \vec{J}
    \quad & \mbox{in $\Omega_1$} \\
  \qquad i \omega\mu \vec{H}^{1,n}+\curl \vec{E}^{1,n} & = & \vec{0}
    & \mbox{in $\Omega_1$} \\
  ({\cal B}_{\vec{n}_1}+{\cal S}_1
  {\cal B}_{\vec{n}_2})(\vec{E}^{1,n},\vec{H}^{1,n})
& = & 
  ({\cal B}_{\vec{n}_1}+{\cal S}_1
  {\cal B}_{\vec{n}_2})(\vec{E}^{2,n-1},\vec{H}^{2,n-1})
    & \mbox{on $\Gamma_{12}$}\\
  - i \omega\varepsilon \vec{E}^{2,n}
      +\curl \vec{H}^{2,n} -\sigma\vec{E}^{2,n} & = & \vec{J}
    & \mbox{in $\Omega_2$} \\
  \qquad i \omega\mu \vec{H}+\curl \vec{E} & = & \vec{0}
    & \mbox{in $\Omega_2$} \\
  ({\cal B}_{\vec{n}_2}+{\cal S}_2
  {\cal B}_{\vec{n}_1})(\vec{E}^{2,n},\vec{H}^{2,n})
    & = & 
  ({\cal B}_{\vec{n}_2}+{\cal S}_2
  {\cal B}_{\vec{n}_1})(\vec{E}^{1,n-1},\vec{H}^{1,n-1})
    & \mbox{on $\Gamma_{21}$,} 
  \end{array}
\end{equation}
where ${\cal S}_j$, $j=1,2$ is a tangential, possibly
pseudo-differential operator we will study in great detail in order to
obtain various optimized Schwarz methods.

\subsection{Convergence Analysis for the Classical Schwarz Algorithm}

We now study properties of the classical Schwarz algorithm
(\ref{ClassicalSchwarz}). We use Fourier analysis, and thus assume
that the coefficients are constant, and the domain on which the
original problem is posed is $\Omega=\R^3$, in which case we need for
Maxwell's equations the Silver-M\"uller radiation condition
\begin{equation} \label{eq1:SM}  
  \lim_{r\rightarrow\infty}r\left(\vec{H}\times \vec{n}-\vec{E}\right)=0,
\end{equation}  
where $r=|\vec{x}|$, $\vec{n}=\vec{x}/{|\vec{x}|}$, in order to obtain
well-posed problems, see \cite{Nedelec:2001:AEE}.  The two subdomains
are now half spaces, 
\begin{equation}\label{TwoHalfPlanes}
  \Omega_1=(0,\infty)\times \R^2,\quad 
  \Omega_2=(-\infty,L)\times \R^2,
\end{equation}
the interfaces are $\Gamma_{12}=\{L\}\times\R^2$ and
$\Gamma_{21}=\{0\}\times\R^2$, and the overlap is $L\ge 0$. We denote
by $k_y$ and $k_z$ the Fourier variables corresponding to a transform
with respect to $y$ and $z$, respectively, and $|\vec{k}|^2 =
k_y^2+k_z^2$.
\begin{theorem}\label{Proposition6}  
  For any given initial guess $(\vec{E}^{1,0};\vec{H}^{1,0})\in
  (L^2(\Omega_1))^6$, $(\vec{E}^{2,0};\vec{H}^{2,0})\in
  (L^2(\Omega_2))^6$, the classical Schwarz algorithm
  (\ref{ClassicalSchwarz}) with overlap $L\ge 0$, including the
  non-overlapping case, is for $\sigma>0$ convergent in
  $(L^2(\Omega_1))^6\times (L^2(\Omega_2))^6$, and the convergence
  factor for each Fourier mode $\vec{k}$ is
  \begin{equation}\label{RTHcla}
    \rho_{cla}(\vec{k},\tilde{\omega},\sigma,Z,L) = \left|
     \frac{\sqrt{|\vec{k}|^2-{\tilde\omega}^2+i\tilde \omega\sigma Z} 
       - i{\tilde\omega}}  
     {\sqrt{|\vec{k}|^2-{\tilde\omega}^2+i\tilde \omega\sigma Z}
       +i{\tilde\omega}}
     e^{-\sqrt{|\vec{k}|^2-{\tilde\omega}^2+i\tilde\omega\sigma Z}L}\right|,
  \end{equation}
  where $\tilde\omega=\omega\sqrt{\varepsilon\mu}$, and $Z =
  \sqrt{\frac{\mu}{\varepsilon}}$ is the impedance as before.
\end{theorem}  
\begin{proof}   
Because of linearity, it suffices to analyze the convergence to the  
zero solution when the right hand side vanishes. Performing a Fourier  
transform of system (\ref{eq:THMaxwell}) in the $y$ and $z$  
direction, the first and the fourth equation provide an algebraic  
expression for $\hat{E}_1$ and $\hat{H}_1$, which is in agreement with  
the fact that these are the characteristic variables associated  
with the null eigenvalue. Inserting these expressions into the remaining  
Fourier transformed equations, we obtain the first order system  
\begin{equation}\label{TH:systemFourier}  
  \arraycolsep0.15em
  \partial_x \left(\begin{array}{c}  
    \hat{E}_2 \\ \hat{E}_3 \\ \hat{H}_2 \\ \hat{H}_3  
  \end{array}\right) +    
  \left[\begin{array}{cccc}  
    0 & 0 & -\frac{k_yk_z}{i\omega\varepsilon+\sigma} & \frac{-{\tilde\omega}^2 +k_y^2+i\omega\mu\sigma}{i\omega\varepsilon+\sigma} \\   
    0 & 0 & \frac{{\tilde\omega}^2-k_z^2-i\omega\mu\sigma}{i\omega\varepsilon+\sigma} & \frac{k_yk_z}{i\omega\varepsilon+\sigma} \\  
    \frac{k_yk_z}{i\omega\mu} & \frac{{\tilde\omega}^2-k_y^2-i\omega\mu\sigma}{i\omega\mu} & 0 & 0 \\  
    \frac{-{\tilde\omega}^2+k_z^2+i\omega\mu\sigma}{i\omega\mu}  & -\frac{k_yk_z}{i\omega\mu} & 0 &0  
  \end{array}\right] \left(\begin{array}{c}  
    \hat{E}_2 \\ \hat{E}_3 \\ \hat{H}_2 \\ \hat{H}_3  
  \end{array}\right) =  
  \left(\begin{array}{c}  
    0 \\ 0 \\0 \\0  
  \end{array}\right).  
\end{equation}  
The eigenvalues of the matrix in (\ref{TH:systemFourier}) and their  
corresponding eigenvectors are  
\begin{equation}\label{v1v2}
  \arraycolsep0.3em
  \lambda_{1,2}^{TH}=-\sqrt{|\vec{k}|^2-{\tilde\omega}^2+i\omega\mu\sigma},
  \quad \vec{v}_1= \left(\begin{array}{c}
  \frac{k_yk_z}{(i{\omega}\varepsilon+\sigma)\lambda}
  \\
  \frac{{-\tilde\omega}^2+k_z^2+i\omega\mu\sigma}{(i{\omega}\varepsilon+\sigma)\lambda}
  \\ 1\\ 0
  \end{array} \right), \quad   
  \vec{v}_2= \left(\begin{array}{c}  
    \frac{{\tilde\omega}^2-k_y^2-i\omega\mu\sigma}{(i{\omega}\varepsilon+\sigma)\lambda} \\   
    -\frac{k_yk_z}{(i{\omega}\varepsilon+\sigma)\lambda} \\   
      0\\  
      1  
  \end{array} \right),  
\end{equation}
and  
\begin{equation}\label{v3v4}
  \arraycolsep0.3em
  \lambda_{3,4}^{TH}=\sqrt{|\vec{k}|^2-{\tilde\omega}^2+i\omega\mu\sigma},\quad  \vec{v}_3=   
  \left(\begin{array}{c}  
    -\frac{k_yk_z}{(i{\omega}\varepsilon+\sigma)\lambda} \\   
    \frac{{\tilde\omega}^2-k_z^2-i\omega\mu\sigma}{(i{\omega}\varepsilon+\sigma)\lambda} \\   
    1\\   
    0  
  \end{array} \right), \quad  
  \vec{v}_4= \left(\begin{array}{c}  
    \frac{k_y^2-{\tilde\omega}^2+i\omega\mu\sigma}{(i{\omega}\varepsilon+\sigma)\lambda} \\   
    \frac{k_yk_z}{(i{\omega}\varepsilon+\sigma)\lambda} \\   
    0\\  
    1  
  \end{array} \right).
\end{equation}
where we set
$\lambda:=\sqrt{|\vec{k}|^2-{\tilde\omega}^2+i\omega\mu\sigma}$.
Because of the radiation condition, the solutions of system
(\ref{TH:systemFourier}) in $\Omega_l$, $l=1,2$, are given by
\begin{equation} \label{SolTH}  
  (\hat E_2^1;\hat E_3^1;\hat H_2^1;\hat H_3^1)=(\alpha_1\vec{v}_1+\alpha_2\vec{v}_2)  
     e^{\lambda(x-L)}, \qquad    
  (\hat E_2^2;\hat E_3^2;\hat H_2^2;\hat H_3^2)=(\beta_1\vec{v}_3+\beta_2\vec{v}_4)   
     e^{-\lambda x},   
\end{equation}  
where the coefficients $\alpha_j$ and $\beta_j$ ($j=1,2$) are uniquely  
determined by the transmission conditions. At the $n$-th step of the  
Schwarz algorithm, the coefficients $\vec{\alpha}=(\alpha_1,\alpha_2)$ and  
$\vec{\beta}=(\beta_1,\beta_2)$ satisfy the system  
$$  
  \vec{\alpha}^{n} = A_1^{-1}A_2e^{-\lambda L}  
    \vec{\beta}^{n-1}, \qquad  
  \vec{\beta}^{n} = B_1^{-1}B_2e^{-\lambda L}  
    \vec{\alpha}^{n-1},  
$$  
where the matrices in the iteration are given by  
\begin{equation} \label{AjTHcl}
\begin{array}{c}  
  A_1=\arraycolsep0.2em
    \left[\begin{array}{cc}  
    -k_yk_z & k_y^2\!-\!{\tilde\omega}^2\!+\!i{\tilde\omega}\lambda+\sigma Z(\lambda+i\tilde\omega) \\   
    k_z^2\!-\!{\tilde\omega}^2\!+\!i{\tilde\omega}\lambda +\sigma Z(\lambda+i\tilde\omega)& -k_yk_z  
  \end{array}\right], \\ 
  A_2 =\left[\begin{array}{cc}  
    k_yk_z & -k_y^2\!+\!{\tilde\omega}^2\!+\!i{\tilde\omega}\lambda+\sigma Z(\lambda-i\tilde\omega) \\   
    -k_z^2\!+\!{\tilde\omega}^2\!+\!i{\tilde\omega}\lambda+\sigma Z(\lambda-i\tilde\omega) & k_yk_z  
  \end{array}\right],  
\end{array}
\end{equation}  
and where $B_l=A_l,l=1,2$.\\
A complete iteration over two steps of the Schwarz algorithm leads
then to
$$  
  \vec{\alpha}^{n+1}=(A_1^{-1}A_2)^2e^{-2\lambda L}\vec{\alpha}^{n-1},  
  \qquad  
  \vec{\beta}^{n+1}= (A_1^{-1}A_2)^2 e^{-2\lambda L}\vec{\beta}^{n-1},  
$$  
and we obtain the iteration matrix 
\begin{equation}\label{ITmat}  
  R = (A_1^{-1}A_2)^2e^{-2\lambda L}=   
  \arraycolsep0.1em
    \left[\begin{array}{cc} 
      \frac{|\vec{k}|^4+2\lambda\sigma Z(k_y^2-k_z^2)+\lambda^2\sigma^2 Z^2}
      {(\lambda+i{\tilde\omega})^2(\lambda+i{\tilde\omega}+\sigma Z)^2} 
      &
      \frac{4k_yk_z\lambda\sigma Z}
      {(\lambda+i{\tilde\omega})^2(\lambda+i{\tilde\omega}+\sigma Z)^2}\\
      \frac{4k_yk_z\lambda\sigma Z}
      {(\lambda+i{\tilde\omega})^2(\lambda+i{\tilde\omega}+\sigma Z)^2} 
      & 
      \frac{|\vec{k}|^4+2\lambda\sigma Z(k_z^2-k_y^2)+\lambda^2\sigma^2 Z^2}
      {(\lambda+i{\tilde\omega})^2(\lambda+i{\tilde\omega}+\sigma Z)^2}
    \end{array}\right] 
    e^{-2\lambda L}. 
\end{equation} 
Now by the definition of $\lambda$, we have
$|\vec{k}|^2=\lambda^2+\tilde\omega^2-i\tilde\omega\sigma Z$, and thus
this matrix can be re-written in factored form,
$$  
  R=\left(\frac{\lambda-i\tilde\omega}{\lambda+i\tilde\omega}\right)^2
  e^{-2\lambda L} Id+\frac{4\lambda\sigma Z}
    {(\lambda + i{\tilde\omega})^2(\lambda + i{\tilde\omega}+\sigma Z)^2}
    \left[\begin{array}{cc}
      -k_z^2 & k_yk_z\\
      k_yk_z & -k_y^2
    \end{array}\right]
  e^{-2\lambda L}  
$$ 
The convergence factor $\rho_{cla}$ of the algorithm is given by the
square root of the spectral radius of the matrix $R$, whose
eigenvalues are
$\left(\frac{\lambda-i\tilde\omega}{\lambda+i\tilde\omega}\right)^2
e^{-2\lambda L}$ and $\left(\frac{\lambda-i\tilde\omega-\sigma
Z}{\lambda+i\tilde\omega+\sigma Z}\right)^2e^{-2\lambda L}$.  Since
$\sigma\ge 0$, a direct computation shows that the convergence factor
is given by the first eigenvalue, which leads to (\ref{RTHcla}), and when
$\sigma \ne 0$, a straightforward computation shows that
$\rho_{cla}(\vec{k}) < 1$ for all Fourier modes $\vec{k}$.
\end{proof}  

If $\sigma=0$, the convergence factor becomes  
\begin{equation} \label{RTH}  
  \rho_{cla}(\vec{k},\tilde{\omega},0,Z,L)=\left\{\begin{array}{ll}  
    \left|\frac{\sqrt{{\tilde\omega}^2-|\vec{k}|^2}-{\tilde\omega}}  
    {\sqrt{{\tilde\omega}^2-|\vec{k}|^2} + {\tilde\omega}}\right|, &   
    \mbox{for $|\vec{k}|^2\le{\tilde\omega}^2$}, \\  
    e^{-\sqrt{|\vec{k}|^2-{\tilde\omega}^2} L},&   
    \mbox{for $|\vec{k}|^2>{\tilde\omega}^2$}.  
  \end{array}\right.  
\end{equation}  
In this case, we obtain for $|\vec{k}|^2={\tilde\omega}^2$ that the
convergence factor equals 1, independently of the overlap, which
indicates that the algorithm is not convergent for $\sigma=0$ when
used in the iterative form described here. In practice, Schwarz
methods are however often used as preconditioners for Krylov methods,
which can handle such isolated problems in the spectrum.  We also see
from the convergence factor (\ref{RTH}) that in the case $\sigma=0$ 
the overlap is necessary for the convergence of the evanescent modes,
$|\vec{k}|^2>{\tilde\omega}^2$. Without overlap, $L=0$, we have
$\rho(|\vec{k}|) <1$ only for the propagative modes,
$|\vec{k}|^2<{\tilde\omega}^2$, and $\rho(|\vec{k}|)=1$ when
$|\vec{k}|^2\geq{\tilde\omega}^2$.

Very similar observations were made in the analysis of optimized Schwarz
methods for the Helmholtz equation in \cite{Gander:2001:OSH}. If one
applies to the Helmholtz equation
\begin{equation}\label{ScalarHelmholtz}
  (\Delta +\tilde{\omega}^2)u=f, \qquad\mbox{in $\Omega=\R^3$},
\end{equation}
with Sommerfeld radiation conditions
$\lim_{r\rightarrow\infty}r\left(\frac{\partial u}{\partial r}
-i\tilde\omega u\right)=0$ and the same two subdomain decomposition
(\ref{TwoHalfPlanes}) the somewhat particular overlapping Schwarz
method (note the unequal treatment in the transmission conditions)
\begin{equation}\label{SchwarzHelm} 
  \arraycolsep0.2em
    \begin{array}{rcllrcll} 
    ({\tilde\omega}^2 + \Delta) u_1^{1,n} & = & f
      \ & \mbox{in $\Omega_1$}  & 
    ({\tilde\omega}^2 + \Delta) u_1^{2,n} & = & f 
      \ &\mbox{in $\Omega_2$},\\ 
    u_1^{1,n} & = & u_1^{2,n-1} & \mbox{on $\Gamma_{12}$}\quad&
    (\partial_x\!-\!i{\tilde\omega}) u_1^{2,n}
    & = & (\partial_x\!-\!i{\tilde\omega}) u_1^{1,n-1}
       & \mbox{on $\Gamma_{21}$},
    \end{array}   
\end{equation}  
then one obtains precisely the same convergence factor (\ref{RTH}).
The classical overlapping Schwarz algorithm with characteristic
transmission conditions (\ref{ClassicalSchwarz}) for Maxwell's
equations is thus very much related to the particular overlapping
Schwarz method (\ref{SchwarzHelm}) for the Helmholtz problem when
$\sigma=0$. This particular Schwarz method is a very simple variant of
an optimized Schwarz method, where one has only replaced one of the
Dirichlet transmission conditions with a better one adapted for low
frequencies. There are much better transmission conditions for
Helmholtz problems, as it was shown in \cite{Gander:2001:OSH}. These
conditions are based on approximations of transparent boundary
conditions, which we will study in the next subsection for Maxwell's
equations.

\subsection{Transparent Boundary Conditions}  \label{transparentSec}
  
To design optimized Schwarz methods for Maxwell's equations, we
derive now transparent boundary conditions for those equations,
following the approach in \cite{Hagstrom:2007:RBC}. We consider the time harmonic
Maxwell's equations (\ref{eq:THMaxwell}) on the domains
$\Omega_1=(-\infty,L)\times \R^2$ and $\Omega_2 =
(0,\infty)\times\R^2$ with right hand sides $\vec{J}_{1,2}$ compactly
supported in $\Omega_{1,2}$, together with the boundary conditions
\begin{equation}\label{ABCmaxw}  
  (\vec{w}^2_+ + {\cal S}_1 \vec{w}^2_-)(0,y,z) = 0, \quad  
  (\vec{w}^1_- + {\cal S}_2 \vec{w}^1_+)(L,y,z) = 0, \qquad  
  (y,z)\in\R^2,  
\end{equation}  
and with Silver-M\"uller condition on their unbounded part, where
$\vec{w}^1_-$ and $\vec{w}^2_+$ are defined in (\ref{W-W+}), and the
operators ${\cal S}_l$, $l=1,2$, are general, pseudo-differential
operators acting in the $y$ and $z$ directions.
  
\begin{theorem}\label{ABCTH}  
If the operators ${\cal S}_l$, $l=1,2$ have the Fourier symbol  
\begin{equation}\label{ABCsymbMTH}  
  \mathcal{F}({\cal S}_l)= \frac{1}{(\lambda+i{\tilde\omega})(\lambda+i{\tilde\omega+\sigma Z})}
    \left[\begin{array}{cc}   
      k_y^2-k_z^2-\lambda\sigma Z & -2k_yk_z \\   
     -2k_yk_z & k_z^2-k_y^2-\lambda\sigma Z  
    \end{array}\right],
\end{equation}  
where $\lambda =\sqrt{|\vec{k}|^2-{\tilde\omega}^2+i\tilde\omega\sigma
Z}$, then the solution of Maxwell's equations (\ref{eq:THMaxwell}) in
$\Omega_{1,2}$ with boundary conditions (\ref{ABCmaxw}) coincides with
the restriction on $\Omega_{1,2}$ of the solution of Maxwell's
equations (\ref{eq:THMaxwell}) on $\R^3$.
\end{theorem}   
\begin{proof}  
We show that the difference $\vec{e}^i,\,i=1,2$ between the solution
of the global problem and the solution of the restricted problem
vanishes. We consider the case of the second domain, similar
computations can be carried out for the first one.  This difference
satisfies in $\Omega_2$ the homogeneous counterpart of
(\ref{eq:THMaxwell}) with homogeneous boundary conditions
(\ref{ABCmaxw}), and we obtain after a Fourier transform in $y$ and
$z$
$$  
  \vec{\hat{e}}^2=(\alpha_1 \vec{v}_1 + \alpha_2 \vec{v}_2)  
    e^{\lambda x}  +   
    (\alpha_3\vec{v}_3+\alpha_4\vec{v}_4)   
    e^{-\lambda x},   
$$ 
where the vectors $\vec{v}_j$, $j=1,..,4$, are defined in (\ref{v1v2})
and (\ref{v3v4}). The Silver-M\"uller radiation condition implies that
$\alpha_1=\alpha_2=0$. Using now the boundary condition
(\ref{ABCmaxw}) at $(0,y,z)$, we obtain that the coefficients
$\alpha_j$, $j=3,4$, satisfy the system of equations
$$  
  (A_1+{\cal S}_1A_2)\left[\begin{array}{c}   
    \alpha_3 \\   
    \alpha_4   
  \end{array}\right] = 0,
$$
where $A_1$ and $A_2$ are defined by (\ref{AjTHcl}). A direct
computation
$$
  \left[\begin{array}{cc}  
    -k_yk_z & k_y^2-{\tilde\omega}^2+i{\tilde\omega}\lambda\\   
    k_z^2-{\tilde\omega}^2+i{\tilde\omega}\lambda &-k_yk_z   
  \end{array} \right]  
  \left[\begin{array}{c}   
    \alpha_3 \\   
    \alpha_4   
  \end{array}\right] =  
  \left[\begin{array}{c}   
    0 \\   
    0   
  \end{array}\right],  
$$  
which implies $\alpha_3=\alpha_4=0$. Thus $\vec{\widehat{e}}^2=\vec{0}$,  
which concludes the proof.  
\end{proof}  
\begin{remark}\label{remark3}
  As in the case of the Cauchy-Riemann equations, see
  \cite{Dolean:2007:WCS}, the symbols in (\ref{ABCsymbMTH}) can be
  written in several, mathematically equivalent forms,
  $$
  \begin{array}{rcccl}
  \mathcal{F}({\cal S}_l)& = & 
  \frac{1}{(\lambda\!+\!i{\tilde\omega})(\lambda+i\tilde\omega+\sigma Z)}M
   & = & \frac{1}{|\vec{k}|^2+\lambda\sigma Z}
    \frac{\lambda\!-\!i\tilde\omega}{\lambda\!+\!i\tilde\omega} M\\
  &=&\frac{1}{|\vec{k}|^2-\lambda\sigma Z}
    \frac{\lambda\!-\!i\tilde\omega-\sigma Z}
    {\lambda\!+\!i\tilde\omega+\sigma Z} M 
  &=&(\lambda\!-\!i{\tilde\omega})(\lambda\!-\!i\tilde\omega-\sigma Z)M^{-1},
  \end{array}
  $$
  where the matrix $M$ is given by 
  $$ 
    M=\left[\begin{array}{cc} 
     k_y^2-k_z^2-\lambda\sigma Z &-2k_yk_z \\ 
     -2k_yk_z & k_z^2-k_y^2-\lambda\sigma Z
    \end{array}\right].
  $$ 
\end{remark}
This motivates different approximations of the transparent conditions
in the context of optimized Schwarz methods. In the case $\sigma = 0$
the first form contains a local and a non-local term, since
multiplication with the matrix $M$ corresponds to second order
derivatives in $y$ and $z$, which are local operations, whereas the
term containing the square-root of $|\vec{k}|^2$ represents a
non-local operation. The last form contains two non-local operations,
since the inversion of the matrix $M$ corresponds to an
integration. This integration can however be passed to the other side
of the transmission conditions by multiplication with the matrix $M$
from the right. The second form contains two non-local terms and a
local one. We propose in the next section several approximations based
on these different forms, and analyze the performance of the
associated optimized Schwarz algorithms.
  
\subsection{Optimized Schwarz Algorithms for Maxwell's Equations}  
  
The transparent operators ${\cal S}_l$, $l=1,2$, introduced in
Subsection \ref{transparentSec}, are important in the development of
optimized Schwarz methods. When used in algorithm
(\ref{OptimizedSchwarz}), they lead to the best possible performance
of the method, as we will show in Remark \ref{RemarkOptimal}.  The
transparent operators are however non-local operators, and hence
difficult to use in practice. In optimized Schwarz methods, they are
therefore approximated to obtain practical methods. If one is willing
to use second order transmission conditions, then the only parts of
the symbols in (\ref{ABCsymbMTH}) that need to be approximated are the
terms $\lambda=\sqrt{|\vec{k}|^2-{\tilde\omega}^2+i\tilde\omega\sigma
Z}$, because the entries of the matrices are polynomials in the
Fourier variables, which correspond to derivatives in the $y$ and $z$
direction.
\begin{theorem}\label{Proposition7}  
  For the optimized Schwarz algorithm (\ref{OptimizedSchwarz}) with
  the two subdomain decomposition (\ref{TwoHalfPlanes}), we obtain for
  $\sigma=0$ the following results:
  \begin{enumerate}
  \item If the operators $\mathcal{S}_1$ and $\mathcal{S}_2$ have the
  Fourier symbol
  \begin{equation}\label{Sigmal1}  
   \sigma_l:=\mathcal{F}(\mathcal{S}_l)= \gamma_l  
   \left[\begin{array}{cc}  
     k_y^2-k_z^2 & -2k_yk_z \\  
     -2k_yk_z & k_z^2-k_y^2  
   \end{array}\right],\quad \gamma_l\in \C(k_z,k_y), \ l=1,2,  
  \end{equation}  
  then the convergence factor is  
  \begin{equation}\label{rhoMTHgen}
  \begin{array}{l}  
    \rho= 
     \textstyle\!\left|\!   
      \frac{(\sqrt{|\vec{k}|^2\!-{\tilde\omega}^2}\!-i{\tilde\omega})^2}  
           {(\sqrt{|\vec{k}|^2\!-{\tilde\omega}^2}\!+i{\tilde\omega})^2}  
      \frac{1\!-\!\gamma_1(\sqrt{|\vec{k}|^2\!-{\tilde\omega}^2}
        \!+i{\tilde\omega})^2}  
           {1\!-\!\gamma_1(\sqrt{|\vec{k}|^2\!-{\tilde\omega}^2}
        \!-i{\tilde\omega})^2}
      \frac{1\!-\!\gamma_2(\sqrt{|\vec{k}|^2\!-{\tilde\omega}^2}
        \!+i{\tilde\omega})^2}  
        {1\!-\!\gamma_2(\sqrt{|\vec{k}|^2\!-{\tilde\omega}^2}
        \!-i{\tilde\omega})^2}
        e^{\!-\!2\sqrt{|\vec{k}|^2\!-\!{\tilde\omega}^2} L}
        \!\right|^{\frac{1}{2}}\!\!\!.  
\end{array}
  \end{equation}  
  \item If the operators $\mathcal{S}_1$ and $\mathcal{S}_2$ have the
  Fourier symbol
  \begin{equation}\label{Sigmal2}  
   \sigma_l:=\mathcal{F}(\mathcal{S}_l)= \delta_l  
   \left[\begin{array}{cc}  
     k_y^2-k_z^2 & -2k_yk_z \\  
     -2k_yk_z & k_z^2-k_y^2  
   \end{array}\right]^{-1},\quad \gamma_l\in \C(k_z,k_y),\ l=1,2,  
  \end{equation}  
  then the convergence factor is  
  \begin{equation}\label{rhoMTHgen2}  
  \begin{array}{l}  
    \rho=
      \!\left|\! 
      \frac{(\sqrt{|\vec{k}|^2\!-{\tilde\omega}^2}+i{\tilde\omega})^2}  
           {(\sqrt{|\vec{k}|^2\!-{\tilde\omega}^2}-i{\tilde\omega})^2}
    \frac{\delta_1-(\sqrt{|\vec{k}|^2\!-{\tilde\omega}^2}-i{\tilde\omega})^2}
         {\delta_1-(\sqrt{|\vec{k}|^2\!-{\tilde\omega}^2}+i{\tilde\omega})^2}
    \frac{\delta_2-(\sqrt{|\vec{k}|^2\!-{\tilde\omega}^2}-i{\tilde\omega})^2}
         {\delta_2-(\sqrt{|\vec{k}|^2\!-{\tilde\omega}^2}+i{\tilde\omega})^2}
      e^{\!-\!2\sqrt{|\vec{k}|^2\!-{\tilde\omega}^2} L}\! 
    \right|^{\frac{1}{2}}\!\!\!.  
\end{array}
  \end{equation}  
  \item If the operator $\mathcal{S}_1$ has the Fourier symbol
  (\ref{Sigmal1}) and $\mathcal{S}_2$ has the Fourier symbol
  (\ref{Sigmal2}), then the convergence factor is
  \begin{equation}\label{rhoMTHgen3}  
    \rho = \textstyle\left|  
    \frac{1-\gamma_1(\sqrt{|\vec{k}|^2-{\tilde\omega}^2}+i{\tilde\omega})^2}  
         {1-\gamma_1(\sqrt{|\vec{k}|^2-{\tilde\omega}^2}-i{\tilde\omega})^2} 
    \frac{\delta_2-(\sqrt{|\vec{k}|^2-{\tilde\omega}^2}-i{\tilde\omega})^2}  
         {\delta_2-(\sqrt{|\vec{k}|^2-{\tilde\omega}^2}+i{\tilde\omega})^2}   
        e^{-2\sqrt{|\vec{k}|^2-{\tilde\omega}^2} L} \right|^{1/2}.  
  \end{equation}  
  \end{enumerate}
\end{theorem}
\begin{proof}  
  The convergence results are again based on Fourier analysis, as in  
  Section \ref{convTHTD}. At the $n$-th step of the Schwarz  
  algorithm, the coefficients $\vec{\alpha}^n=(\alpha^{1,n},\alpha^{2,n})$  
  and $\vec{\beta}=(\beta_1,\beta_2)$ in (\ref{SolTH}) satisfy  
  \begin{equation}\label{IterSchw}  
    \vec{\alpha}^{n}=\bar{A}_1^{-1}\bar{A}_2  
      e^{-\lambda L}\vec{\beta}^{n-1}, \qquad  
    \vec{\beta}^{n}=\bar{B}_1^{-1}\bar{B}_2\ e^{-\lambda L}\vec{\alpha}^{n-1},  
  \end{equation}  
  where $\lambda=\sqrt{|\vec{k}|^2-{\tilde\omega}^2}$, and
  the matrices $\bar{A}_l$ and $\bar{B}_l$, $l=1,2$, are given by  
  $$  
    \bar{A}_1 = A_1+\sigma_1 A_2, \quad  
    \bar{A}_2 = A_2+\sigma_1 A_1,\quad  
    \bar{B}_1 = A_1+\sigma_2 A_2, \quad  
    \bar{B}_2 = A_2+\sigma_2 A_1,
  $$  
  with $A_l$, $l=1,2$, defined in (\ref{AjTHcl}). A complete double
  iteration of the Schwarz algorithm leads therefore to
  $$  
    \vec{\alpha}^{n+1} = \bar{A}_1^{-1}\bar{A}_2\bar{B}_1^{-1}\bar{B}_2  
      e^{-2\lambda L}\vec{\alpha}^{n-1}, \qquad  
    \vec{\beta}^{n+1} = \bar{B}_1^{-1}\bar{B}_2\bar{A}_1^{-1}\bar{A}_2  
      e^{-2\lambda L}\vec{\beta}^{n-1}.  
  $$ 
  Notice that the matrices $A_1$ and $A_2$ verify the properties
  \begin{equation}\label{prop}
    MA_1=-(\lambda+i\tilde\omega)^2A_2,
    \quad MA_2=-(\lambda-i\tilde\omega)^2A_1, 
  \end{equation}
  which are essential in all three cases:
  \begin{enumerate}
  \item We obtain, using (\ref{prop}),
  $$
    \begin{array}{rcccl}
      \bar A_1&=&A_1+\gamma_1MA_2
        &=&(1-(\lambda-i\tilde\omega)^2\gamma_1)A_1,\\
      \bar A_2&=&A_2+\gamma_1MA_1
        &=&(1-(\lambda+i\tilde\omega)^2\gamma_1)A_2,\\
      \bar B_1&=&A_1+\gamma_2MA_2
        &=&(1-(\lambda-i\tilde\omega)^2\gamma_2)A_1,\\
      \bar B_2&=&A_2+\gamma_2MA_1
        &=&(1-(\lambda+i\tilde\omega)^2\gamma_2)A_2,
    \end{array}
  $$
  and therefore the iteration matrix becomes  
  $$ 
    R_1 = \bar{A}_1^{-1}\bar{A}_2\bar{B}_1^{-1}\bar{B}_2
    = \textstyle \frac{(1-\gamma_1(\lambda+i{\tilde\omega})^2)
      (1-\gamma_2(\lambda+i{\tilde\omega})^2)}
      {(1- \gamma_1(\lambda-i{\tilde\omega})^2)
      (1-\gamma_2(\lambda-i{\tilde\omega})^2)}
      (A_1^{-1}A_2)^2 e^{-2\lambda L},
  $$ 
  and by using the spectral radius of (\ref{ITmat}) the result
  follows.  
  \item We get
  $$
    \begin{array}{rcccl}
      \bar A_1 &=& A_1+\delta_1M^{-1}A_2 
       &=&\left(1-\frac{\delta_1}{(\lambda+i\tilde\omega)^2}\right)A_1,\\
      \bar A_2 &=& A_2+\delta_1M^{-1}A_1
       &=& \left(1-\frac{\delta_1}{(\lambda-i\tilde\omega)^2}\right)A_2,\\
      \bar B_1 &=& A_1+\delta_2M^{-1}A_2 
       &=& \left(1-\frac{\delta_2}{(\lambda+i\tilde\omega)^2}\right)A_1,\\
      \bar B_2 &=& A_2+\delta_2M^{-1}A_1 
       &=& \left(1-\frac{\delta_2}{(\lambda-i\tilde\omega)^2}\right)A_2,
\end{array}
$$
and thus the iteration matrix is
  $$  
    R_2=\bar{A}_1^{-1}\bar{A}_2\bar{B}_1^{-1}\bar{B}_2 
    = \textstyle \left(\frac{\lambda+i{\tilde\omega}}
    {\lambda-i{\tilde\omega}}\right)^4   
    \frac{(\delta_1-(\lambda-i{\tilde\omega})^2)
    (\delta_2-(\lambda-i{\tilde\omega})^2)}  
    {(\delta_1 - (\lambda+i{\tilde\omega})^2)
    (\delta_2-(\lambda+i{\tilde\omega})^2)}
    (A_1^{-1}A_2)^2 e^{-2\lambda L},  
  $$  
  and we use again the spectral radius of (\ref{ITmat}) to conclude.
  \item  The conclusion follows as in the first two cases.
  \end{enumerate}
\end{proof}  
\begin{remark}\label{RemarkOptimal}  
  From (\ref{rhoMTHgen}), we see that the choice
  $\gamma_1=\gamma_2=1/(\sqrt{|\vec{k}|^2-{\tilde\omega}^2}
  +i{\tilde\omega})^2$ is optimal, since then $\rho\equiv 0$, for
  all frequencies $\vec{k}$. With this choice of $\gamma_1$ and
  $\gamma_2$, the matrices $\bar{A}_2$ and $\bar{B}_2$ actually
  vanish.
\end{remark}  
  
\subsection{A hierarchy of optimized transmission conditions}
\label{CasesSection}  
  
We present now several particular choices of the remaining parameters
in the transmission operators $\mathcal{S}_l$ in Theorem
\ref{Proposition7}. To facilitate the use of our results in domain
decomposition codes, we return to the initial notation using the
physical parameters $\omega$, $\varepsilon$ and $\mu$.
\begin{description}  
  
\item[Case 1:] taking $\gamma_1=\gamma_2=0$ in (\ref{Sigmal1}), which
  amounts to enforce the classical characteristic Dirichlet
  transmission conditions, the convergence factor is
  $$ 
    \rho_{1}({\omega,\varepsilon,\mu},L,|\vec{k}|) = \left|\left(
      \frac{\sqrt{|\vec{k}|^2-\omega^2\varepsilon\mu}-i{\omega\sqrt{\varepsilon\mu}}}
      {\sqrt{|\vec{k}|^2-\omega^2\varepsilon\mu}+i{\omega\sqrt{\varepsilon\mu}}}\right)^2
      e^{-2\sqrt{|\vec{k}|^2-\omega^2\varepsilon\mu} L} \right|^{\frac{1}{2}}.
  $$  
  In the non-overlapping case, $L=0$, this choice ensures convergence
  only for propagative modes, and corresponds to the Taylor
  transmission conditions of order zero proposed in the seminal paper
  \cite{Despres:1993:DDM} for the Helmholtz equation.
  
\item[Case 2:] taking
  $\gamma_1=\gamma_2=\frac{1}{|\vec{k}|^2}
  \frac{s-i{\omega\sqrt{\varepsilon\mu}}}{s+i{\omega\sqrt{\varepsilon\mu}}}$ in (\ref{Sigmal1}) or
  $\gamma_1=\frac{1}{|\vec{k}|^2-2\omega^2\varepsilon\mu +2 i{\omega\sqrt{\varepsilon\mu}} s}$
  in (\ref{Sigmal1}) and $\delta_2={|\vec{k}|^2-2\omega^2\varepsilon\mu
  -2 i{\omega\sqrt{\varepsilon\mu}} s}$ in (\ref{Sigmal2}) with $s\in\mathbb{C}$, the
  convergence factor is
  $$  
    \rho_{2}({\omega,\varepsilon,\mu},L,|\vec{k}|,s) = \left|\left(
      \frac{\sqrt{|\vec{k}|^2-\omega^2\varepsilon\mu}-s}{\sqrt{|\vec{k}|^2-\omega^2\varepsilon\mu}+s}
      \right)^2 e^{-2\sqrt{|\vec{k}|^2-\omega^2\varepsilon\mu} L} \right|^{\frac{1}{2}}.  
  $$  
  which is for $L=0$ identical to the convergence factor obtained for
  optimized non-overlapping Schwarz methods for the Helmholtz equation
  in \cite{Gander:2001:OSH}.

\item[Case 3:] taking
  $\gamma_1=\gamma_2=\frac{1}{|\vec{k}|^2-2\omega^2\varepsilon\mu +2
  i{\omega\sqrt{\varepsilon\mu}} s}$ in (\ref{Sigmal1}) with $s\in\mathbb{C}$, the
  convergence factor is
  $$  
  \begin{array}{lcl}
    \rho_{3}({\omega,\varepsilon,\mu},L,|\vec{k}|,s) & = & \left|
     \frac{\sqrt{|\vec{k}|^2-\omega^2\varepsilon\mu}-i{\omega\sqrt{\varepsilon\mu}}}
       {\sqrt{|\vec{k}|^2-\omega^2\varepsilon\mu}+i{\omega\sqrt{\varepsilon\mu}}}\right|
      \rho_{2}({\omega,\varepsilon,\mu},L,|\vec{k}|,s) \\[2ex]
     & \le & \rho_{2}(\omega,\varepsilon,\mu,L,|\vec{k}|,s).
  \end{array} 
  $$

\item[Case 4:] taking
  $\gamma_{l}=\frac{1}{|\vec{k}|^2}\frac{s_{l}-i{\omega\sqrt{\varepsilon\mu}}}{s_l+i{\omega\sqrt{\varepsilon\mu}}},\,l=1,2$
  in (\ref{Sigmal1}) or
  $\gamma_1=\frac{1}{|\vec{k}|^2-2\omega^2\varepsilon\mu +2
  i{\omega\sqrt{\varepsilon\mu}} s_1}$ in (\ref{Sigmal1}) and
  $\delta_2={|\vec{k}|^2-2\omega^2\varepsilon\mu -2
  i{\omega\sqrt{\varepsilon\mu}} s_2}$ in (\ref{Sigmal2}) with
  $s_l\in\mathbb{C},\,l=1,2$, the convergence factor is
  $$  
    \rho_{4}({\omega,\varepsilon,\mu},L,|\vec{k}|,s_1,s_2) = \textstyle\left|
    \frac{\sqrt{|\vec{k}|^2-\omega^2\varepsilon\mu}-s_1}{\sqrt{|\vec{k}|^2-\omega^2\varepsilon\mu}+s_1}
    \frac{\sqrt{|\vec{k}|^2-\omega^2\varepsilon\mu}-s_2}{\sqrt{|\vec{k}|^2-\omega^2\varepsilon\mu}+s_2}   
      e^{-2\sqrt{|\vec{k}|^2-\omega^2\varepsilon\mu} L} \right|^{\frac{1}{2}},  
  $$  
  which is for $L=0$ identical to the convergence factor obtained for
  a two sided non-overlapping optimized Schwarz method for the
  Helmholtz equation in \cite{gander-etal:06}.

\item[Case 5:] taking $\gamma_l=\frac{1}{|\vec{k}|^2-2\omega^2\varepsilon\mu +2
  i{\omega\sqrt{\varepsilon\mu}} s_l}$ in (\ref{Sigmal1}) with $s_l\in\mathbb{C},\,l=1,2$,
  the convergence factor is
  $$  
  \begin{array}{lcl}
    \rho_{5}({\omega,\varepsilon,\mu},L,|\vec{k}|,s_1,s_2) & = & \left|
     \frac{\sqrt{|\vec{k}|^2-\omega^2\varepsilon\mu}-i{\omega\sqrt{\varepsilon\mu}}}
        {\sqrt{|\vec{k}|^2-\omega^2\varepsilon\mu}+i{\omega\sqrt{\varepsilon\mu}}}\right|
        \rho_{4}({\omega,\varepsilon,\mu},L,|\vec{k}|,s_1,s_2) \\[2ex]
      & \le &  
        \rho_{4}({\omega,\varepsilon,\mu},L,|\vec{k}|,s_1,s_2). 
\end{array}
  $$
\end{description} 
Except for Case 1, all cases use second order transmission conditions,
even though we use only a zeroth order approximation of the non-local
operator $\sqrt{|\vec{k}|^2-\omega^2\varepsilon\mu}$. 
In the cases with parameters, the best choice for the parameters is in
general the one that minimizes the convergence factor for all
$|\vec{k}|\in K$, where $K$ denotes the set of relevant numerical
frequencies.  
One therefore needs to solve the min-max problems
\begin{equation}  \label{minmaxproblems}
  \min_{s\in\C}\max_{|\vec{k}|\in K}
     \rho_j({\omega,\varepsilon,\mu},L,|\vec{k}|,s),\ j=2,3,\quad
  \min_{s_1,s_2\in C}\max_{|\vec{k}|\in K}
     \rho_j({\omega,\varepsilon,\mu},L,|\vec{k}|,s_1,s_2),\ j=4,5.  
\end{equation}
We can choose $K = [(k_{\min},k_-)\cup (k_+,k_{\max})]^2$, where
$k_{min}$ denotes the smallest frequency relevant to the subdomain,
and $k_{max}=\frac{C}{h}$ denotes the largest frequency supported by the
numerical grid with mesh size $h$, and $k_{\pm}$ are parameters to be
chosen to exclude the resonance frequencies. If for example the
domain $\Omega$ is a rectilinear conductor with homogeneous Dirichlet
conditions on the lateral surface, the solution is the sum of the
transverse electric (TE) and transverse magnetic (TM) fields. If the
transverse section of the conductor is a rectangle with sides of
length $a$ and $b$, the TE and TM fields can be expanded in a Fourier
series with the harmonics $\sin(\frac{n\pi y}{a})\sin(\frac{m\pi
z}{b})$, where the relevant frequencies are
$|\vec{k}|=\pi\sqrt{\frac{m^2}{a^2}+\frac{n^2}{b^2}}$, $m,n\in \N^+$.
The lowest one is therefore
$k_{min}=\pi\sqrt{\frac{1}{a^2}+\frac{1}{b^2}}$, and if the mesh size
$h$ satisfies $h=\frac{a}{N}=\frac{b}{M}$, where $N$ and $M$ are the
number of grid points in the $y$ and $z$ direction, then the highest
frequency would be $k_{\max}=\frac{\sqrt{2}\pi}{h}$. The parameters
$k_{\pm}$ would correspond to the frequencies closest to ${\omega\sqrt{\varepsilon\mu}}$,
i.e.  $k_-=\pi\sqrt{\frac{m_1^2}{a^2}+\frac{n_1^2}{b^2}}$ and
$k_+=\pi\sqrt{\frac{m_2^2}{a^2}+\frac{n_2^2}{b^2}}$, where
$\pi\sqrt{\frac{m_1^2}{a^2}+\frac{n_1^2}{b^2}}<{\omega\sqrt{\varepsilon\mu}}<
\pi\sqrt{\frac{m_2^2}{a^2}+\frac{n_2^2}{b^2}}$, but such precise
estimates are not necessary if Krylov acceleration is used, see
\cite{Gander:2001:OSH,gander-etal:06}.

The complete mathematical analysis of the min-max problems
(\ref{minmaxproblems}) is hard, and currently open for $L>0$. When
$L=0$, i.e. no overlap, Case 2 and Case 4 are equivalent to the
corresponding optimized Schwarz method for the Helmholtz equation, for
which theoretical results are available, see
\cite{gander-etal:06}. Here, we use asymptotic analysis and an
equioscillation principle to solve all the min-max problems in
(\ref{minmaxproblems}) asymptotically as the mesh size goes to zero,
in order to obtain compact formulas for the best parameters to be used
in our numerical simulations. This leads to the asymptotic formulas
for the optimized parameters of the form $s=p(1-i)$ and
$s_l=p_l(1-i)$, $l=1,2$, with $p$ and $p_l$ shown in Table
\ref{TabOptMax}.
\begin{table}  
  \centering  
  \tabcolsep0.1em
  \begin{tabular}{|c|c|c|c|c|}\hline  
            & \multicolumn{2}{c|}{with overlap, $L=h$} & 
              \multicolumn{2}{c|}{without overlap, $L=0$} \\ \hline   
     Case   & $\rho$ & parameters & $\rho$ & parameters  \\ \hline   
      1 & $1-\sqrt{k_+^2-\tilde\omega^2}h$ & 
       none
       & $1$ & 
       none\\  
      2 & $1-2C_{\tilde\omega}^{\frac{1}{6}}h^{\frac{1}{3}}$ & 
       $p = \frac{C_{\tilde\omega}^{\frac{1}{3}}}{2\cdot h^{\frac{1}{3}}}$ &
       $1-\frac{\sqrt{2}C_{\tilde\omega}^{\frac{1}{4}}}{\sqrt{C}}\sqrt{h}$ & 
       $p=\frac{\sqrt{C}C_{\tilde\omega}^{\frac{1}{4}}}{\sqrt{2}\sqrt{h}}$\\  
      3 & $1-2(k_+^2-\tilde\omega^2)^{\frac{1}{6}}h^{\frac{1}{3}}$ & 
       $p = \frac{(k_+^2-\tilde\omega^2)^{\frac{1}{3}}}{2\cdot h^{\frac{1}{3}}}$
       & $1\!-\!\frac{\sqrt{2}(k_+^2-\tilde\omega^2)^{\frac{1}{4}}}{\sqrt{C}}\sqrt{h}$ & 
       $ p=\frac{\sqrt{C}(k_+^2-\tilde\omega^2)^{\frac{1}{4}}}{\sqrt{2}\sqrt{h}}$\\  
      4 & $1-2^{\frac{2}{5}}C_{\tilde\omega}^{\frac{1}{10}}h^{\frac{1}{5}}$ &
      $ 
      \left\{\begin{array}{l}
         p_1 = \frac{C_{\tilde\omega}^{\frac{2}{5}}}{2^{\frac{7}{5}}\cdot h^{\frac{1}{5}}},\\
         p_2 = \frac{C_{\tilde\omega}^{\frac{1}{5}}}{2^{\frac{6}{5}}\cdot h^{\frac{3}{5}}}
      \end{array}\right.
      $
      &$1-\frac{C_{\tilde\omega}^{\frac{1}{8}}}{C^{\frac{1}{4}}}h^{\frac{1}{4}}$
      & $\left\{\begin{array}{l} p_1 =
      \frac{C_{\tilde\omega}^{\frac{3}{8}}\cdot C^{\frac{1}{4}}}{2\cdot
      h^{\frac{1}{4}}},\\ 
      p_2 = \frac{C_{\tilde\omega}^{\frac{1}{8}}\cdot
      C^{\frac{3}{4}}}{h^{\frac{3}{4}}}\end{array}\right.$\\ 
      5 & $1\!-\!
      2^{\frac{2}{5}}(k_+^2\!-\!\tilde\omega^2)^{\frac{1}{10}}h^{\frac{1}{5}}$ &
      $\left\{\arraycolsep0.1em\begin{array}{l} p_1 =
      \frac{(k_+^2-\tilde\omega^2)^{\frac{2}{5}}}{2^{\frac{7}{5}}\cdot
      h^{\frac{1}{5}}},\\ p_2 =
      \frac{(k_+^2-\tilde\omega^2)^{\frac{1}{5}}}{2^{\frac{6}{5}}\cdot
      h^{\frac{3}{5}}}\end{array}\right.$ & $1-
      \frac{(k_+^2-\tilde\omega^2)^{\frac{1}{8}}}{C^{\frac{1}{4}}}h^{\frac{1}{4}}$
      & $\left\{\arraycolsep0.1em\begin{array}{l} p_1 =
      \frac{(k_+^2-\tilde\omega^2)^{\frac{3}{8}}\cdot
      C^{\frac{1}{4}}}{2\cdot h^{\frac{1}{4}}},\\ p_2 =
      \frac{(k_+^2-\tilde\omega^2)^{\frac{1}{8}}\cdot
      C^{\frac{3}{4}}}{h^{\frac{3}{4}}}\end{array}\right.$ \\ \hline
  \end{tabular}  
  \caption{Asymptotic convergence factor and optimal choice of the
    parameters in the transmission conditions for the five variants of
    the optimized Schwarz method applied to Maxwell's equations, when
    the mesh parameter $h$ is small, and the maximum numerical
    frequency is estimated by $k_{\max}=\frac{C}{h}$. Here
    $\tilde\omega = \omega\sqrt{\varepsilon\mu}$ and
    $C_{\tilde\omega}=\min\left(k_+^2-\tilde\omega^2,
    \tilde\omega^2-k_-^2\right)$.}
  \label{TabOptMax}  
\end{table}  
These results allow us to compare the performance of all the optimized
Schwarz methods for Maxwell's equations theoretically: we obtain a
hierarchy of better and better convergence factors starting with Case
1 and ending with Case 5. In addition, the explicit formulas for the
optimized parameters can be used in order to easily obtain black-box
optimized Schwarz methods for Maxwell's equations, which would not be
possible otherwise. In Section \ref{NumericalSection}, we will
furthermore verify these theoretical results numerically.

\section{The Case of Time Discretization}\label{TimeDiscretizedSection}
  
If we do not assume the wave to be periodic in time, the time domain
also needs to be discretized. We consider a uniform time grid with
time step $\Delta t$, and use an implicit time integration scheme
for the time derivative in (\ref{maxw}) of the form
$$  
\left\{\begin{array}{l}
  -\sigma \left( \frac{\vec{E}^{n+1}+\vec{E}^{n}}{2} \right)-\varepsilon\frac{\vec{E}^{n+1}-\vec{E}^{n}}{\Delta t}   
    + \curl \left( \frac{\vec{H}^{n+1}+\vec{H}^{n}}{2} \right) = \vec{J},\\ 
  \mu\frac{\vec{H}^{n+1}-\vec{H}^{n}}{\Delta t}   
    + \curl \left( \frac{\vec{E}^{n+1}+\vec{E}^{n}}{2} \right) = \vec{0},  
\end{array}\right.
$$ 
where the mean value is introduced to ensure energy conservation, see
\cite{Dolean:2005:IFV}. With this time discretization, we have to
solve at each time step the system
\begin{equation} \label{eq:TDmaxwsyst}
  -\sigma \vec{E}-\varepsilon\sqrt{\eta}\vec{E}   
    + \curl \vec{H} = \tilde{\vec{J}},\quad 
  \mu\sqrt{\eta}\vec{H}   
    + \curl \vec{E} = \vec{g},  
\end{equation}  
where we have set $(\vec{E},\vec{H}):=(\vec{E}^{n+1},\vec{H}^{n+1})$,
$\sqrt{\eta}:=\frac{2}{\Delta t}$, $\vec{\widetilde{J}}:=\vec{J} -
\sqrt{\eta}\varepsilon\vec{E}^n+2\sigma\vec{E}^n -\curl \vec{H}^n$,
and $\vec{g}=\sqrt{\eta}\mu\vec{H}^n -\curl \vec{E}^n$. 

\subsection{Classical and Optimized Schwarz Algorithm}
  \label{convTD}  

 As in the time
harmonic case, we consider the problem (\ref{eq:TDmaxwsyst}) in a
bounded domain $\Omega$, with either Dirichlet conditions on the
tangent electric field, or impedance conditions, on $\partial \Omega$,
in order to obtain a well posed problem, see \cite{Nedelec:2001:AEE}.
For the two subdomain decomposition in Figure \ref{domaindec}, the
classical Schwarz algorithm would at each time step then perform the
iteration
\begin{equation}\label{SchwarzMaxwTD}  
  \begin{array}{rcll}  
    -\sqrt{\eta}\varepsilon\vec{E}^{1,n}
      +\curl\vec{H}^{1,n}-\sigma\vec{E}^{1,n} &=& \vec{J}^1
    \quad & \mbox{in $\Omega_1$} \\
  \qquad \sqrt{\eta}\mu \vec{H}^{1,n}+\curl \vec{E}^{1,n} & = & \vec{g}^1
    & \mbox{in $\Omega_1$} \\
  {\cal B}_{\vec{n}_1}(\vec{E}^{1,n},\vec{H}^{1,n})& = & 
  {\cal B}_{\vec{n}_1}(\vec{E}^{2,n-1},\vec{H}^{2,n-1})
    & \mbox{on $\Gamma_{12}$}\\
  - \sqrt{\eta}\varepsilon \vec{E}^{2,n}
      +\curl \vec{H}^{2,n} -\sigma\vec{E}^{2,n} & = & \vec{J}^2
    & \mbox{in $\Omega_2$} \\
  \qquad \sqrt{\eta}\mu \vec{H}+\curl \vec{E} & = & \vec{g}^2
    & \mbox{in $\Omega_2$} \\
  {\cal B}_{\vec{n}_2}(\vec{E}^{2,n},\vec{H}^{2,n})
    & = & 
  {\cal B}_{\vec{n}_2}(\vec{E}^{1,n-1},\vec{H}^{1,n-1})
    & \mbox{on $\Gamma_{21}$}. 
  \end{array}
\end{equation}  
\begin{theorem}  
  Let $\Omega=\R^3$ be decomposed into $\Omega_1:=(-\infty, L)\times
  \R^2$ and $\Omega_2:=(0,+\infty) \times \R^2$, $L\geq 0$. Then, for
  any initial guess $(\vec{E}^{1,0};\vec{H}^{1,0})\in
  (L^2(\Omega_1))^6$, $(\vec{E}^{2,0};\vec{H}^{2,0})\in
  (L^2(\Omega_2))^6$, the classical Schwarz algorithm
  (\ref{SchwarzMaxwTD}) with overlap $L\ge 0$, including the
  non-overlapping case, is for $\sigma\ge 0$ convergent in
  $(L^2(\Omega_1))^6\times (L^2(\Omega_2))^6$, and the convergence
  factor is bounded by 
  \begin{equation} \label{RTD}  
     R_{cla}=\frac{\sqrt{L \tilde\eta +2\sqrt{\varepsilon\mu}} -
       \sqrt{L\tilde\eta}} {\sqrt{L \tilde\eta
       +2\sqrt{\varepsilon\mu}} + \sqrt{L\tilde\eta}}e^{-\sqrt{L\eta}
       \sqrt{L\tilde\eta+2\sqrt{\varepsilon\mu}}}<1,
  \end{equation}  
  where $\tilde{\eta}=\eta\varepsilon\mu$. 
\end{theorem}  
\begin{proof}   
  This result follows like in the time harmonic case, simply replacing
  $i{\omega}$ by $\sqrt{\eta}$. The convergence factor in Fourier is 
  $$  
    \rho_{cla}(|\vec{k}|)=\left|\frac{\sqrt{|\vec{k}|^2
      +\eta\varepsilon\mu+\sqrt{\eta}\sigma Z} - \eta\sqrt{\varepsilon\mu}}  
      {\sqrt{|\vec{k}|^2+\eta\varepsilon\mu
      +\sqrt{\eta}\sigma Z}+\eta\sqrt{\varepsilon\mu}}
      e^{-\sqrt{|\vec{k}|^2+\eta\varepsilon\mu+\sqrt{\eta}\sigma Z}L}
      \right|,  
  $$  
  and the method thus converges for all Fourier modes. To conclude the
  proof, it suffices to take the maximum of the convergence factor
  over $|\vec{k}|$.
\end{proof}  

The preceding theorem shows that the classical Schwarz algorithm with
Dirichlet transmission conditions applied to the time-discretized
Maxwell's equations is convergent for all frequencies $|\vec{k}|$, and
that the overlap is not necessary to ensure convergence. The classical
Schwarz algorithm corresponds in the case $\sigma=0$ to a simple
optimized Schwarz algorithm for an associated positive definite
Helmholtz equation
\begin{equation}\label{posdefHelmholtz}
    (\tilde\eta -\Delta) u = f,
\end{equation}
and from \cite{Gander:2006:OSM} we know that there are much
better transmission conditions for such problems. We thus propose at
each time step the new algorithm
\begin{equation}\label{SchwarzOpt}  
  \begin{array}{rcll}  
    -\sqrt{\eta}\varepsilon\vec{E}^{1,n}
      +\curl\vec{H}^{1,n}-\sigma\vec{E}^{1,n} &=& \vec{J}^1
    \quad & \mbox{in $\Omega_1$} \\
  \qquad \sqrt{\eta}\mu \vec{H}^{1,n}+\curl \vec{E}^{1,n} & = & \vec{g}^1
    & \mbox{in $\Omega_1$} \\
  ({\cal B}_{\vec{n}_1}+{\cal S}_1
  {\cal B}_{\vec{n}_2})(\vec{E}^{1,n},\vec{H}^{1,n})& = & 
  ({\cal B}_{\vec{n}_1}+{\cal S}_1
  {\cal B}_{\vec{n}_2})(\vec{E}^{2,n-1},\vec{H}^{2,n-1})
    & \mbox{on $\Gamma_{12}$}\\
  - \sqrt{\eta}\varepsilon \vec{E}^{2,n}
      +\curl \vec{H}^{2,n} -\sigma\vec{E}^{2,n} & = & \vec{J}^2
    & \mbox{in $\Omega_2$} \\
  \qquad \sqrt{\eta}\mu \vec{H}+\curl \vec{E} & = & \vec{g}^2
    & \mbox{in $\Omega_2$} \\
  ({\cal B}_{\vec{n}_2}+{\cal S}_2
  {\cal B}_{\vec{n}_1})(\vec{E}^{2,n},\vec{H}^{2,n})
    & = & 
  ({\cal B}_{\vec{n}_2}+{\cal S}_2
  {\cal B}_{\vec{n}_1})(\vec{E}^{1,n-1},\vec{H}^{1,n-1})
    & \mbox{on $\Gamma_{21}$}. 
  \end{array}
\end{equation}  
Now Theorem \ref{ABCTH}, Remark \ref{remark3}, Theorem
\ref{Proposition7} and all the cases in subsection \ref{CasesSection}
hold unchanged for the time discretized case of Maxwell's equations
upon replacing $i\omega$ by $\sqrt{\eta}$, so we do not restate these
results here. The nature of the associated min-max problems
(\ref{minmaxproblems}) however changes fundamentally, and the
optimization parameters are now real, $s=p\in \R$ and $s_l=p_l\in \R$,
$l=1,2$. For cases 2 and 4 a complete analysis is available, see
\cite{Gander:2006:OSM}. Using a lengthy asymptotic analysis again, we complete
the results for the other cases, and show in Table \ref{TabOptTime}
the asymptotically optimal parameters to use in the time domain case.
\begin{table}  
  \centering  
  \tabcolsep0.1em
  \begin{tabular}{|c|c|c|c|c|}\hline  
            & \multicolumn{2}{c|}{with overlap, $L=h$} & 
              \multicolumn{2}{c|}{without overlap, $L=0$} \\ \hline   
     Case   & $\rho$ & parameters & $\rho$ & parameters  \\ \hline   
      1 & $1-2^{\frac{3}{2}}{\tilde\eta}^{\frac{1}{4}}\sqrt{h}$ & 
       none
       & $1-2\frac{\sqrt{{\tilde\eta}}}{C}h$ & 
       none\\  
      2 & $1-2^{\frac{13}{6}}{\tilde\eta}^{\frac{1}{6}}h^{\frac{1}{3}}$ & 
       $p = \frac{2^{-\frac{1}{3}}{\tilde\eta}^{\frac{1}{3}}}{h^{\frac{1}{3}}}$ &
       $1-\frac{4{\tilde\eta}^{\frac{1}{4}}\sqrt{h}}{\sqrt{C}}$ & 
       $p=\frac{\sqrt{C}{\tilde\eta}^{\frac{1}{4}}}{\sqrt{h}}$\\  
      3 & $1-2^{\frac{7}{4}}{\tilde\eta}^{\frac{1}{8}}h^{\frac{1}{4}}$ & 
       $p = \frac{\sqrt{2}{\tilde\eta}^{\frac{1}{4}}}{\sqrt{h}}$
       & $1-\frac{2^{\frac{5}{3}}{\tilde\eta}^{\frac{1}{6}}}{C^{\frac{1}{3}}}
        h^{\frac{1}{3}}$ & 
       $p=\frac{2^{\frac{2}{3}}C^{\frac{2}{3}}{\tilde\eta}^{\frac{1}{6}}}{h^{\frac{2}{3}}}$\\  
      4 & $1-2^{\frac{4}{5}}{\tilde\eta}^{\frac{1}{10}}
        h^{\frac{1}{5}}$ & 
      $p_1=\frac{{\tilde\eta}^{\frac{1}{5}}}
           {2^{\frac{2}{5}}h^{\frac{3}{5}}},
      p_2=\frac{{\tilde\eta}^{\frac{2}{5}}}{16^{\frac{1}{5}}
       h^{\frac{1}{5}}}$
     &$1-\frac{\sqrt{2}{\tilde\eta}^{\frac{1}{8}}}{C^{\frac{1}{4}}}
      h^{\frac{1}{4}}$ &
    $p_1=\frac{\sqrt{2}C^{\frac{3}{4}}{\tilde\eta}^{\frac{1}{8}}}
     {h^{\frac{3}{4}}},
     p_2=\frac{C^{\frac{1}{4}}{\tilde\eta}^{\frac{3}{8}}}
     {\sqrt{2}h^{\frac{1}{4}}}$\\  
      5 & $1-2^{\frac{7}{6}}{\tilde\eta}^{\frac{1}{12}}h^\frac{1}{6}$  & 
    $p_1=\frac{2^{\frac{2}{3}}{\tilde\eta}^{\frac{1}{3}}}
     {h^{\frac{1}{3}}}$, 
    $p_2=\frac{2^{\frac{1}{3}}{\tilde\eta}^{\frac{1}{6}}}
     {h^{\frac{2}{3}}}$
     &$1 - \frac{2{\tilde\eta}^{\frac{1}{10}}}{C^{\frac{1}{5}}}
      h^\frac{1}{5}$ &
    $p_1=\frac{2C^{\frac{4}{5}}{\tilde\eta}^{\frac{1}{10}}}
     {h^{\frac{4}{5}}}$,
    $p_2=\frac{2C^{\frac{2}{5}}{\tilde\eta}^{\frac{3}{10}}}
     {h^{\frac{2}{5}}}$    
     \\  \hline
  \end{tabular}  
  \caption{Asymptotic convergence factor and optimal choice of the
  parameters in the transmission conditions for the five variants of
  the optimized Schwarz method applied to the time domain Maxwell's
  equations, when the mesh parameter $h$ is small, and the maximum
  numerical frequency is estimated by $k_{\max}=\frac{C}{h}$. Here $\tilde\eta = \eta\varepsilon\mu$.}
  \label{TabOptTime}  
\end{table}
Again we obtain an entire hierarchy of optimized Schwarz methods, with
better and better convergence factors from Case 1 up to Case 5. While
for the time harmonic equations Case 2 and 3, and Case 4 and 5 were
asymptotically comparable, here all cases are asymptotically
different.  It is also interesting to note a relationship of the
optimized parameters for the time domain case with the one for the
Cauchy-Riemann equations, see \cite{Dolean:2007:WCS}: Case 2 and 4 are
identical, since the corresponding convergence rates in the two cases
are the same, while for Case 1, 3 and 5 there is a small difference in
the constants, which is due to the additional low frequency term in
the Maxwell case. The difference appears to be systematic, the
convergence factor of the Maxwell case is obtained from the
convergence factor of the Cauchy-Riemann case by replacing $h$ by
$2h$, while for the optimized parameters one has to multiply by $2$ in
addition to the replacement of $h$ by $2h$.

\section{Numerical Experiments}\label{NumericalSection}

We discretize the equations using a finite volume method on a
staggered grid, which leads to the Yee scheme in the interior. For the
first two test cases we consider the propagation in vacuum with
$\varepsilon=\mu = 1$ and $\sigma=0$. We
first show the two dimensional problem of transverse electric waves,
since this allows us to compute with finer mesh sizes and thus to
illustrate our asymptotic results by numerical experiments.  We
simulate directly the error equations, $f=0$, on a uniform mesh with
mesh parameter $h$, and we use a random initial guess to ensure that
all the frequency components are present in the iteration.  We then
show the full 3d case, first for a model problem, and then for
the application of heating a chicken in a microwave oven.

\subsection{Two-dimensional case}

We consider the transverse electric waves problem (TE) in the plane
$(x, y,0)$. There is no more dependence on $z$ and the components
$E_3$, $H_1$ $H_2$ are identically zero. The problem obtained is formally
identical to the three-dimensional case (\ref{maxwsyst}), if ${\vec u}
= (E_1, E_2, H_3)^t$, and the matrix $N_{\bf v}$ becomes 
  \begin{equation*}
    N_{\bf v} =
    \begin{pmatrix}
      -v_y \\
      v_x
    \end{pmatrix},
  \end{equation*}
and the matrices $G_x$, $G_y$ and  $G_{\bf v}$ are
  \begin{equation*}
    G_x =
    \begin{pmatrix}
       & N_{{\bf e}_x} \\
      N_{{\bf e}_x}^t & 
    \end{pmatrix}
    ,\quad
    G_y =
    \begin{pmatrix}
       & N_{{\bf e}_y} \\
      N_{{\bf e}_y}^t & 
    \end{pmatrix}\quad \text{and} \quad 
    G_{\bf v} = 
    \begin{pmatrix}
       & N_ {\bf v}\\
      N_{\bf v}^t & 
    \end{pmatrix}.
  \end{equation*}
All the analytical results remain valid, we only need to replace
$|\vec{k}|$ by $|k_y|$, and the corresponding quantities in the
optimized parameters for both time-harmonic and time-discretized
solutions. We solve Maxwell's equations on the unit square
$\Omega=(0,1)^2$, decomposed into the two subdomains
$\Omega_1=(0,\beta)\times (0,1)$ and $\Omega_2=(\alpha,1)\times
(0,1)$, where $0<\alpha\le\beta<1$, and therefore the overlap is
$L=\beta-\alpha$, and we consider both decompositions with and without
overlap.

In the time-harmonic case, the frequency ${\tilde\omega}=2\pi$ is chosen such
that the rule of thumb of $10$ points per wavelength is not
violated. Table \ref{TabCompM} shows the iteration count for all
Schwarz algorithms we considered, in the overlapping and
non-overlapping case. 
\begin{table} 
  \centering  
  \begin{tabular}{|c|c|c|c|c|c|c|c|c|}\hline  
            & \multicolumn{4}{c|}{with overlap, $L=h$} & 
              \multicolumn{4}{c|}{without overlap, $L=0$} \\ \hline   
      h & $1/16$ & $1/32$ & $1/64$ & $1/128$  
      &  $1/16$ & $1/32$ & $1/64$ & $1/128$  \\ \hline   
      Case 1&18(17)&27(21)&46(27)&71(33)& -(48)& -(73)& -(100)& -(138)\\ 
      Case 2&16(13)&16(14)&17(15)&20(17)&28(22)&36(26)&50(34)&68(40)\\
      Case 3&10(12)&12(13)&14(14)&16(17)&31(20)&40(23)&56(25)&81(28)\\
      Case 4&17(13)&17(14)&20(16)&22(18)&26(20)&28(24)&33(28)&38(30)\\
      Case 5&10(12)&12(13)&14(15)&17(18)&41(24)&53(26)&63(30)&73(32)\\ \hline
  \end{tabular}  
  \caption{Number of iterations in the 2d time harmonic case to attain
  an error tolerance of $=10^{-6}$ for different transmission
  conditions and different mesh sizes. }
  \label{TabCompM}  
\end{table}
The results are presented in the form $it_{S}(it_{GM})$, where $it_S$
denotes the iteration number for the iterative version of the
algorithm and $it_{GM}$ the iteration number for the accelerated
version using GMRES.

\begin{figure}
\begin{center}
\begin{tabular}{cc}
\includegraphics[scale=0.35]{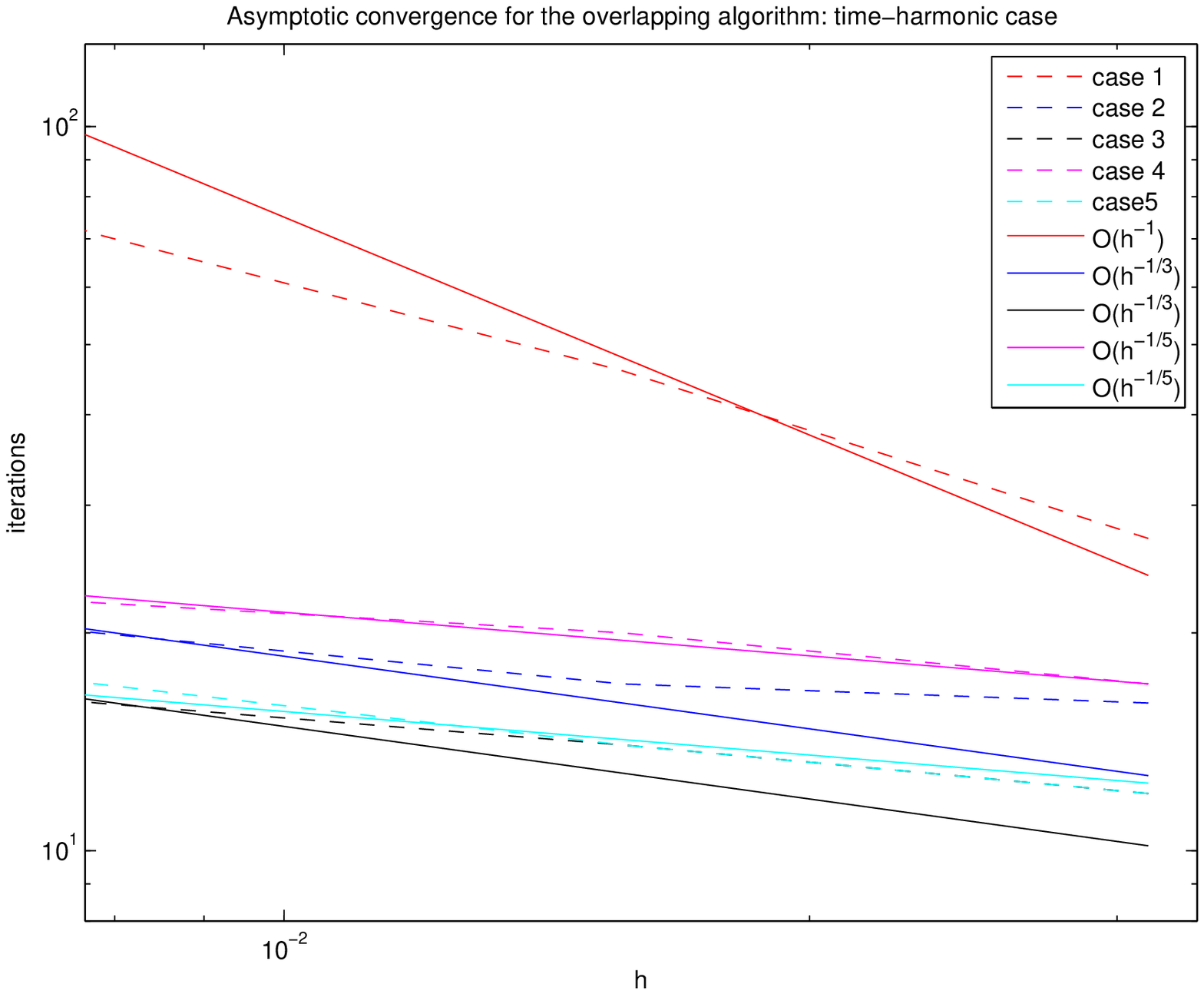} &
\includegraphics[scale=0.35]{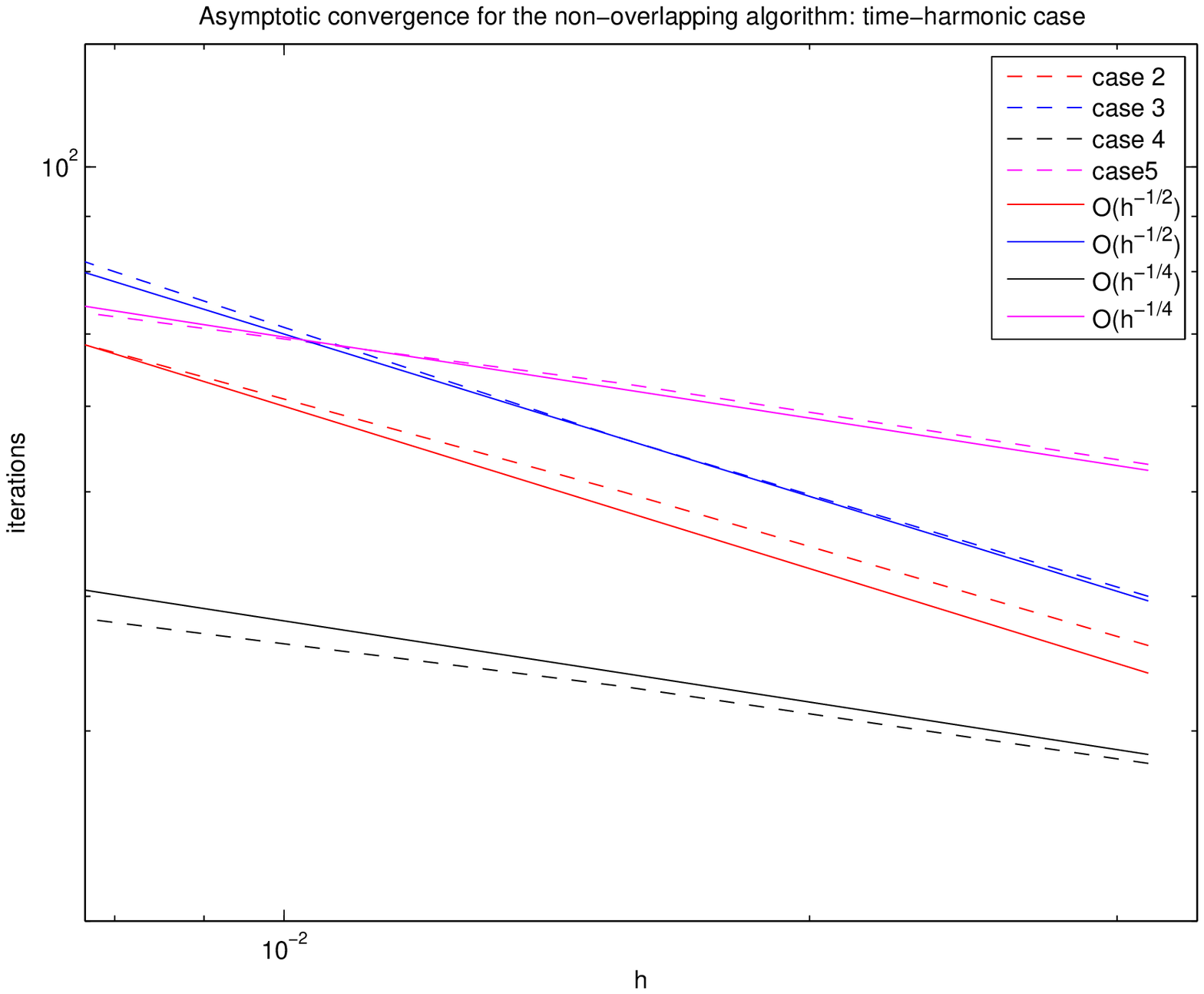} 
\end{tabular}
\end{center}
\caption{Asymptotics for the overlapping (left) and non-overlapping
(right) cases for the time harmonic equations.}
\label{FigHarm}
\end{figure}

In Figure \ref{FigHarm} we show the results we obtained in a graph,
together with the expected asymptotics. 
Both on the left in the overlapping case and on the right in the
non-overlapping one, the asymptotics agree quite well, except for the
classical case with overlap, where the algorithm performs better than
predicted by the asymptotic analysis. In the case of the
Cauchy-Riemann equations \cite{Dolean:2007:WCS}, it was observed
that certain discretizations of the hyperbolic system can introduce
higher order terms in the discretized transmission conditions, which
can improve the convergence behavior, as we observe it here, an 
issue that merits further study.
 
For the time discretized Maxwell's equations we choose ${\tilde\eta}=1$. Table
\ref{TabCompMT} shows the iteration count for all Schwarz algorithms
we considered, in the overlapping and non-overlapping case. We
observe that the classical non-overlapping algorithm converges only
very slowly, the need of optimized methods is evident here.
\begin{table} 
  \centering  
  \begin{tabular}{|c|c|c|c|c|c|c|c|c|}\hline  
            & \multicolumn{4}{c|}{with overlap, $L=h$} & 
              \multicolumn{4}{c|}{without overlap, $L=0$} \\ \hline   
      h & $1/16$ & $1/32$ & $1/64$ & $1/128$  
      &  $1/16$ & $1/32$ & $1/64$ & $1/128$  \\ \hline   
      Case 1 & 17 & 24 & 33 & 45 & 280 & 559 & 1310 & 2630  \\ 
      Case 2 & 13 & 15 & 19 & 24 & 39 & 56 & 77 & 111 \\
      Case 3 & 12 & 14 & 16 & 18 & 13 & 16 & 20 & 26 \\
      Case 4 & 12 & 13 & 15 & 17 & 21 & 25 & 30 & 36 \\
      Case 5 & 12 & 14 & 16 & 18 & 13 & 17 & 19 & 22  
     \\  \hline
  \end{tabular}  
  \caption{Number of iterations in the 2d time discretized case to
  attain an error tolerance of $10^{-6}$ for different transmission
  conditions and different mesh sizes.}
  \label{TabCompMT}  
\end{table}  

\begin{figure}[ht]
\begin{center}
\begin{tabular}{cc}
\includegraphics[scale=0.35]{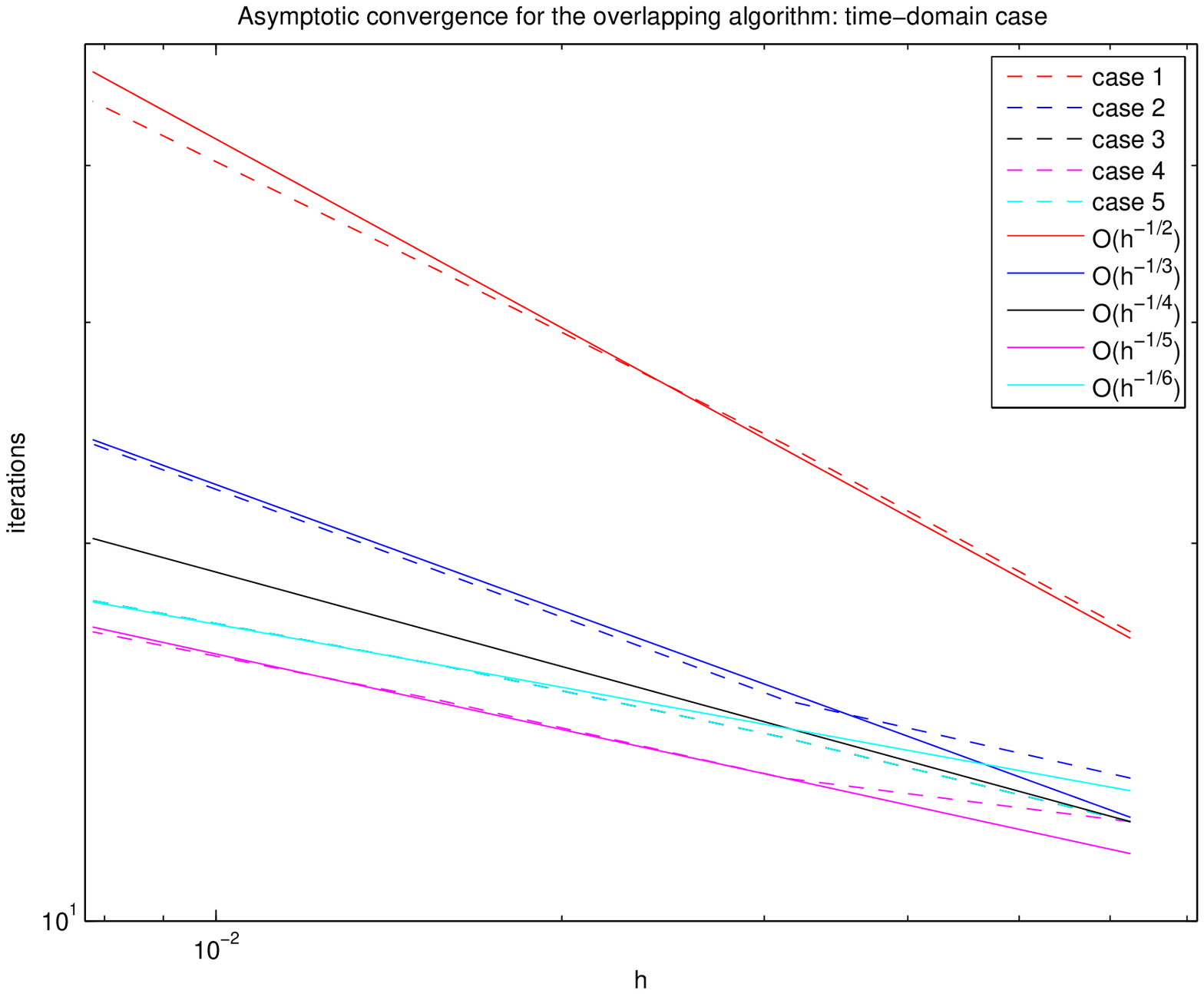} &
\includegraphics[scale=0.35]{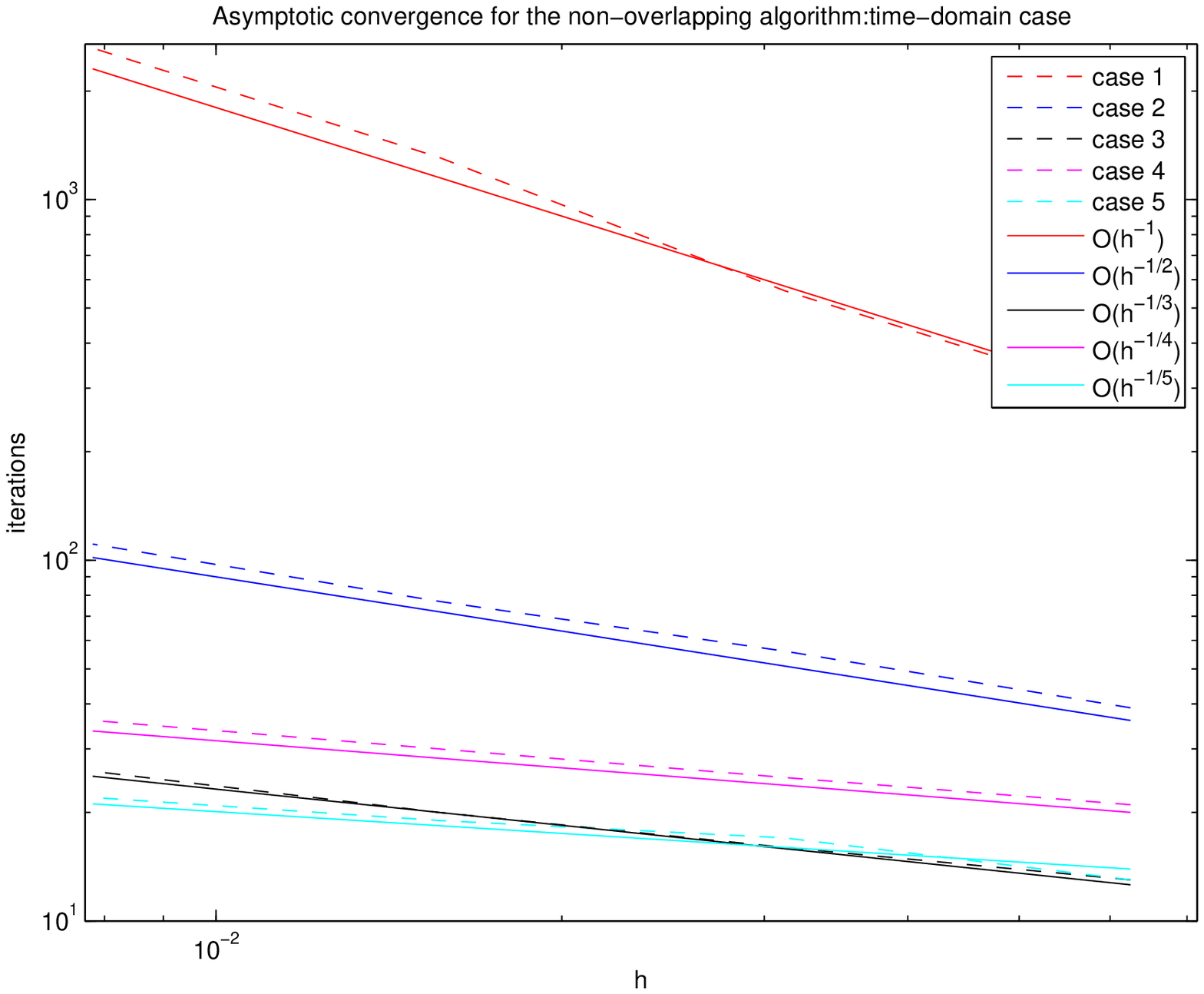} 
\end{tabular}
\end{center}
\caption{Asymptotics for the overlapping (left) and non-overlapping
(right) cases for the time discretized equations.}
\label{FigTime}
\end{figure}

In Figure \ref{FigTime} we show the results we obtained in a graph,
together with the expected asymptotics, and there is very good agreement.
 
\subsection{Three-dimensional case}

We solve now Maxwell's equations on the unit cube $\Omega=(0,1)^3$. We
decompose the domain into two subdomains $\Omega_1=(0,\beta)\times
(0,1)^2$ and $\Omega_2=(\alpha,1)\times (0,1)^2$, with
$0<\alpha\le\beta<1$, and $L=\beta-\alpha$ as before.
In the time-harmonic case, we chose the frequency ${\tilde\omega}=2\pi/3$
to satisfy the rule of thumb of 10 points per wavelength.  Table
\ref{TabCompM3} shows the iteration count for all Schwarz algorithms
we considered, both in the overlapping and non-overlapping
case. 
\begin{table}
  \centering  
  \begin{tabular}{|c|c|c|c|c|c|c|}\hline  
            & \multicolumn{3}{c|}{with overlap, $L=h$} & 
              \multicolumn{3}{c|}{without overlap, $L=0$} \\ \hline   
      h & $1/8$ & $1/16$ & $1/30$   
      &  $1/8$ & $1/16$ & $1/30$   \\ \hline   
      Case 1 & 19(13) & 29(17) & 46(22) & -(93)  & -(140) & -(202) \\ 
      Case 2 & 14(12) & 19(14) & 23(16) & 48(29) & 69(36) & 98(48)  \\
      Case 3 & 16(12) & 18(14) & 21(16) & 65(35) & 80(42) & 166(55)  \\
      Case 4 & 15(13) & 19(15) & 22(17) & 38(28) & 60(33) & 104(39)  \\
      Case 5 & 16(13) & 18(14) & 21(16) & 70(36) & 80(42) & 176(55)  \\  \hline
  \end{tabular}  
  \caption{Number of iterations in the 3d time harmonic case to attain
  an error level of $10^{-6}$ for different transmission conditions
  and different mesh sizes.}
  \label{TabCompM3}  
\end{table}

The results for the time discretized Maxwell's equations where
${\tilde\eta}=1$ are shown in Table \ref{TabCompMT3}.
\begin{table}
  \centering  
  \begin{tabular}{|c|c|c|c|c|c|c|}\hline  
            & \multicolumn{3}{c|}{with overlap, $L=h$} & 
              \multicolumn{3}{c|}{without overlap, $L=0$} \\ \hline   
      h & $1/8$ & $1/16$ & $1/30$   
      &  $1/8$ & $1/16$ & $1/30$  \\ \hline   
      Case 1 & 14 & 18 & 25 & 246 & 467 & 859 \\ 
      Case 2 & 13 & 18 & 22 & 46 & 65 & 87\\
      Case 3 & 12 & 15 & 17 & 47 & 59 & 73 \\
      Case 4 & 14 & 17 & 19 & 48 & 57 & 66 \\
      Case 5 & 12 & 14 & 16 & 46 & 53 & 60   
     \\  \hline
  \end{tabular}  
  \caption{Number of iterations in the 3d time discretized case to
  attain an error level of $10^{-6}$ for different transmission
  conditions and different mesh sizes.}
  \label{TabCompMT3}  
\end{table}  

\subsection{A three-dimensional application: chicken in a micro-wave
  oven}

\begin{figure}
  \centering
  \includegraphics[width=0.49\textwidth]{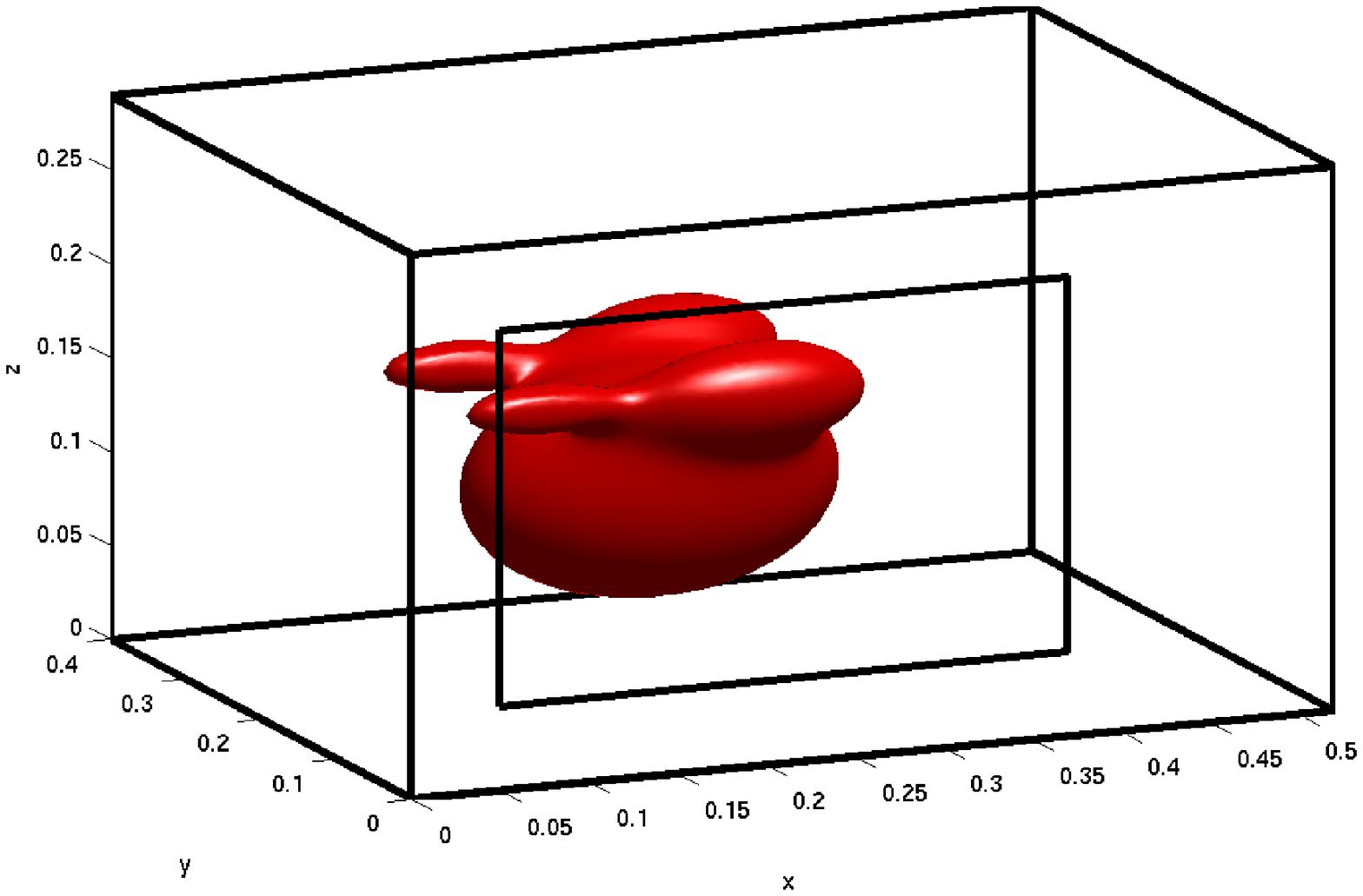}
  \includegraphics[width=0.49\textwidth]{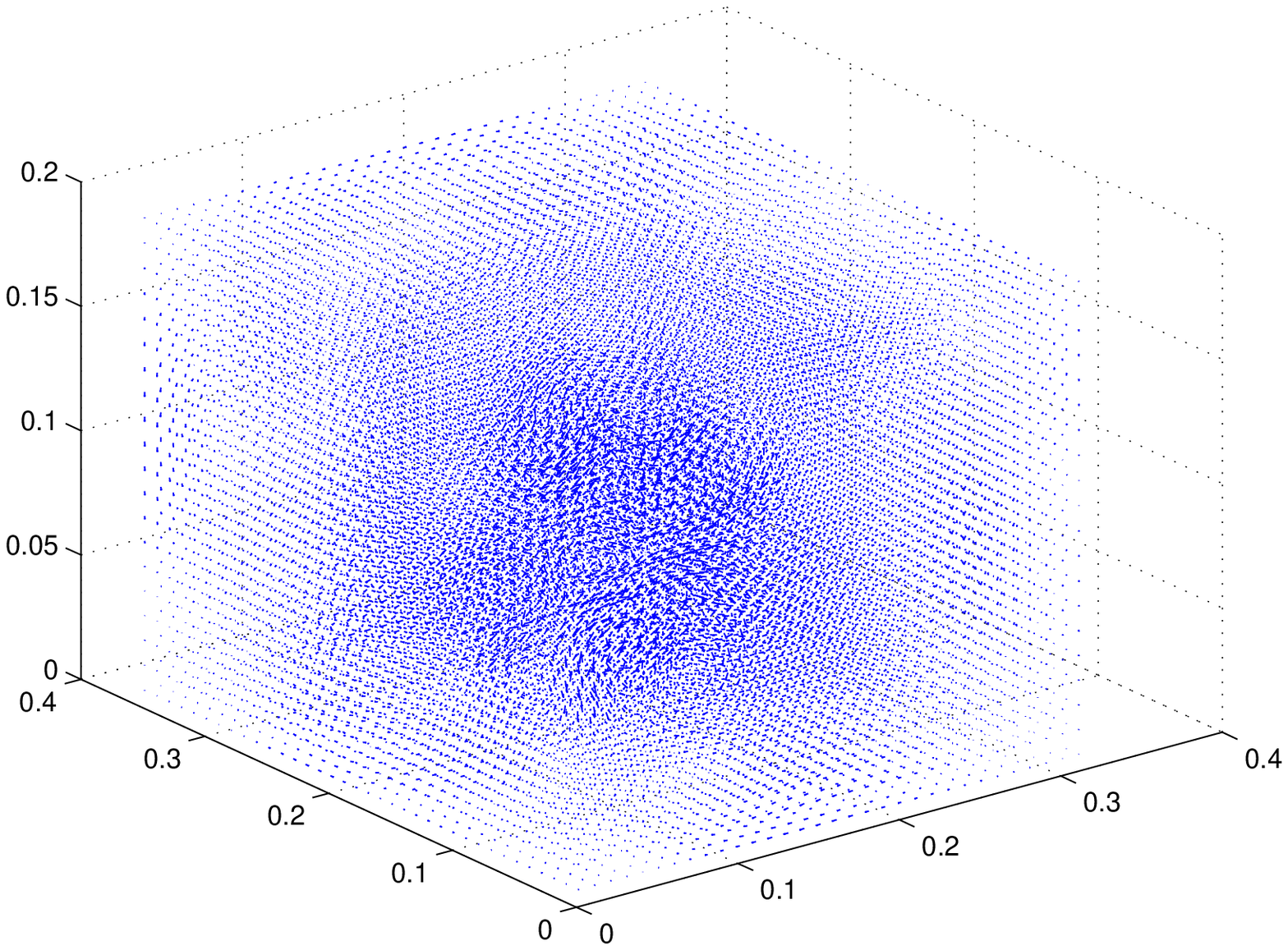}
  \caption{Chicken in a Whirlpool Talent Combi 4 micro-wave oven on
  the left, and real part of the magnetic field in the cooking cavity
  while heating the chicken on the right.}
  \label{Chicken}
\end{figure}

We apply now the previous principles to derive an efficient
domain-decomposition method based on optimized interface conditions to
solve a realistic application: heating up a chicken in a micro-wave
oven, see Figure \ref{Chicken} on the left. The computational domain
is now given by the heating cavity of a Whirlpool Talent Combi 4
microwave oven, $\Omega = [0,0.32] \times [0,0.36] \times [0,0.20]$
meters.
We impose metallic boundary conditions (which means a null tangential
electric field) on all faces except on the right of the oven, where
the components of the electric field are the dominant TE10 mode
generated by the magnetron on a small rectangle of dimensions
$0.08\times 0.04$. The electric and electromagnetic properties of the
media are now non-constant in the computational domain: inside the
chicken, we have an electric permittivity $\varepsilon = 4.43\cdot
10^{-11}\frac{\mbox{\tiny Farads}}{\mbox{\tiny m}}$ and the
conductivity is $\sigma = 3\cdot 10^{-11}\frac{\mbox{\tiny
Siemens}}{\mbox{\tiny m}}$, whereas for the air $\varepsilon =
8,85\cdot 10^{-12} \frac{\mbox{\tiny Farads}}{\mbox{\tiny m}}$ and
$\sigma = 0\frac{\mbox{\tiny Siemens}}{\mbox{\tiny m}}$. The magnetic
permeability is the same for both, $\mu = 4\pi \cdot 10^{-7}
\frac{\mbox{\tiny Henry}}{\mbox{\tiny m}}$, and the frequency is given
by $\omega = 2\pi\cdot 2.45\mbox{ GHz}$.

\begin{figure}
  \centering
  \includegraphics[width=0.49\textwidth]{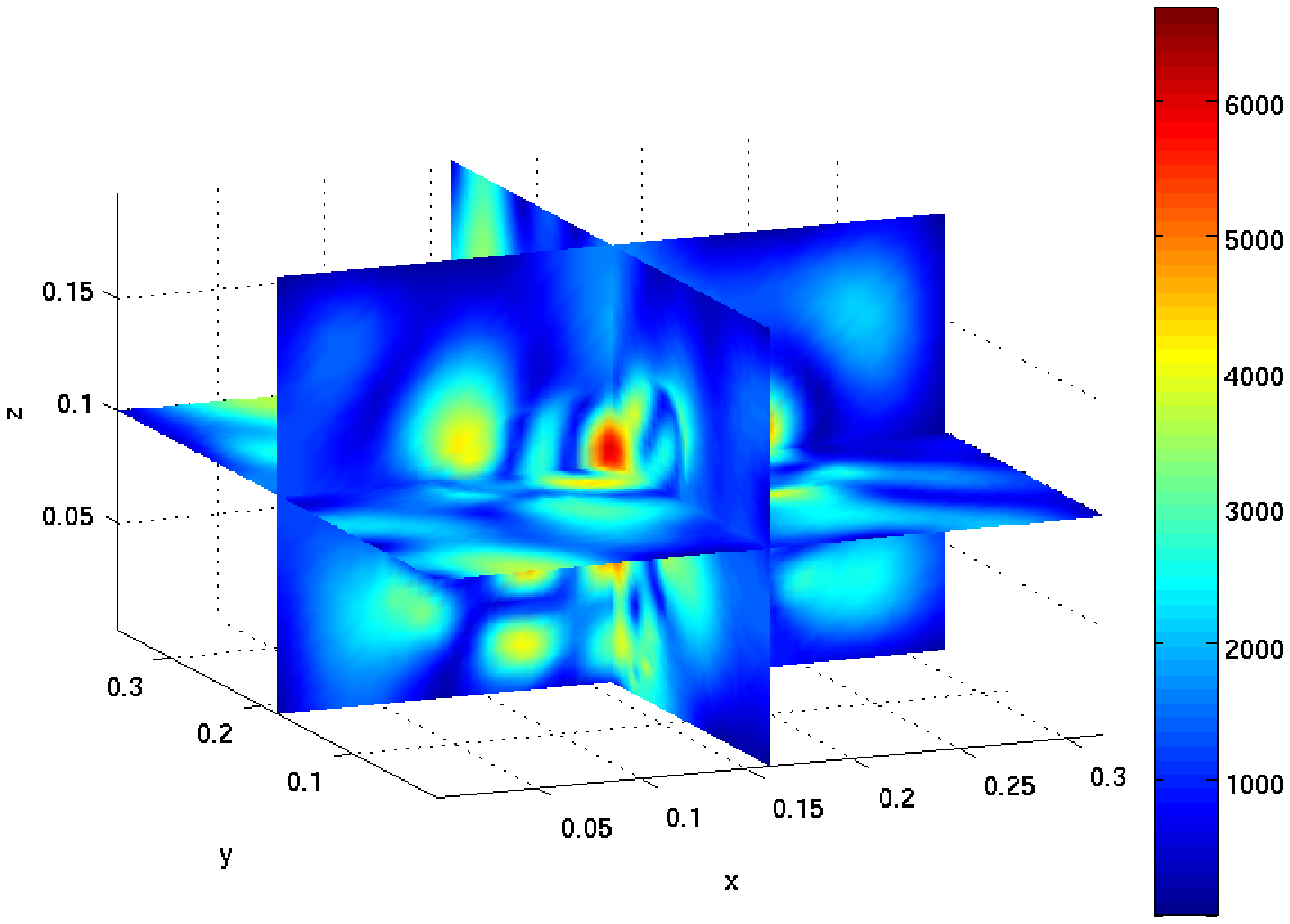}
  \includegraphics[width=0.49\textwidth]{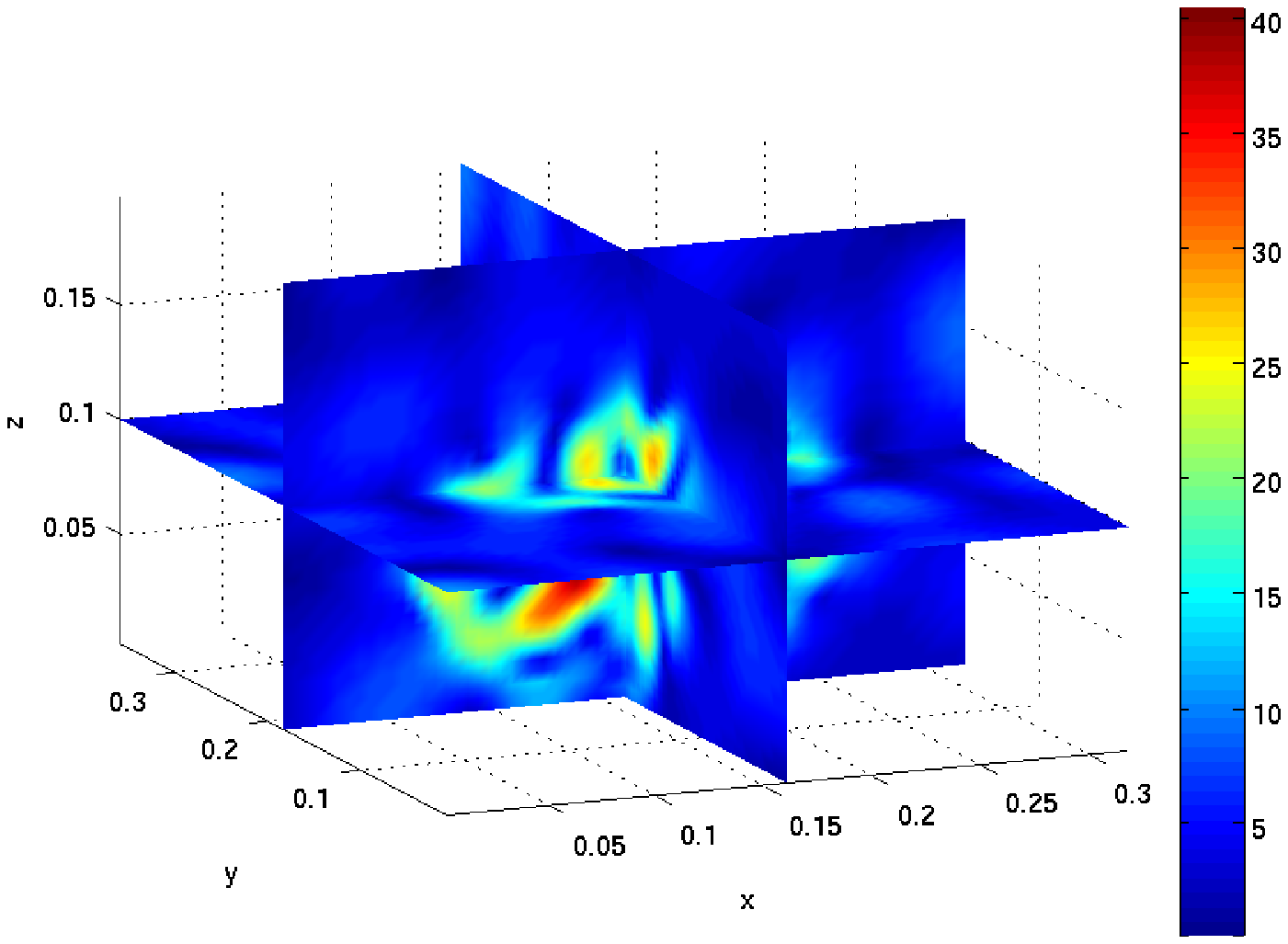} 
  \caption{Chicken heating in a microwave oven: electric field
  intensity on the left, and magnetic field intensity on the right.}
  \label{sol3d}
\end{figure}

\begin{figure}
  \centering
  \includegraphics[width=0.49\textwidth]{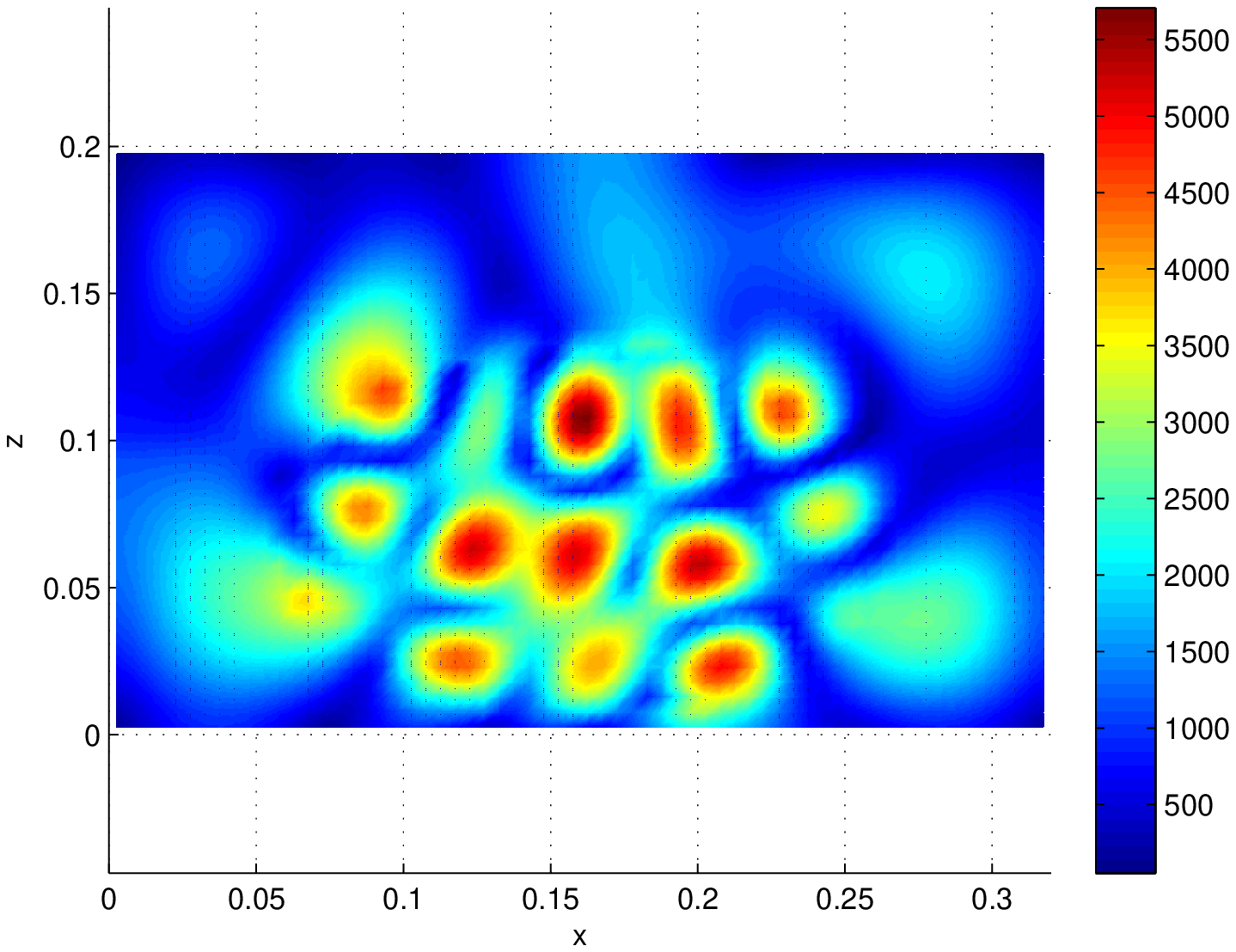} 
  \includegraphics[width=0.49\textwidth]{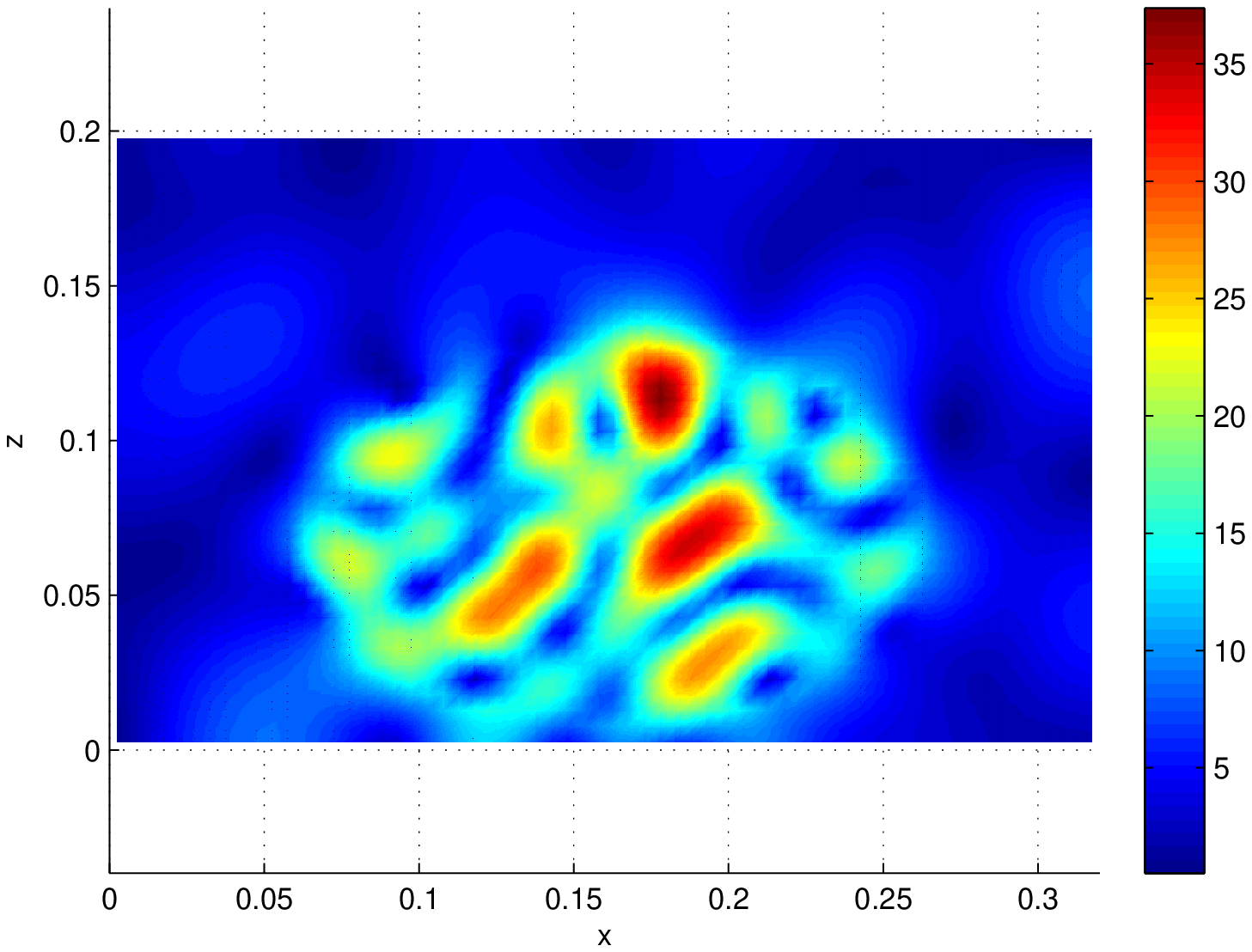}
  \includegraphics[width=0.49\textwidth]{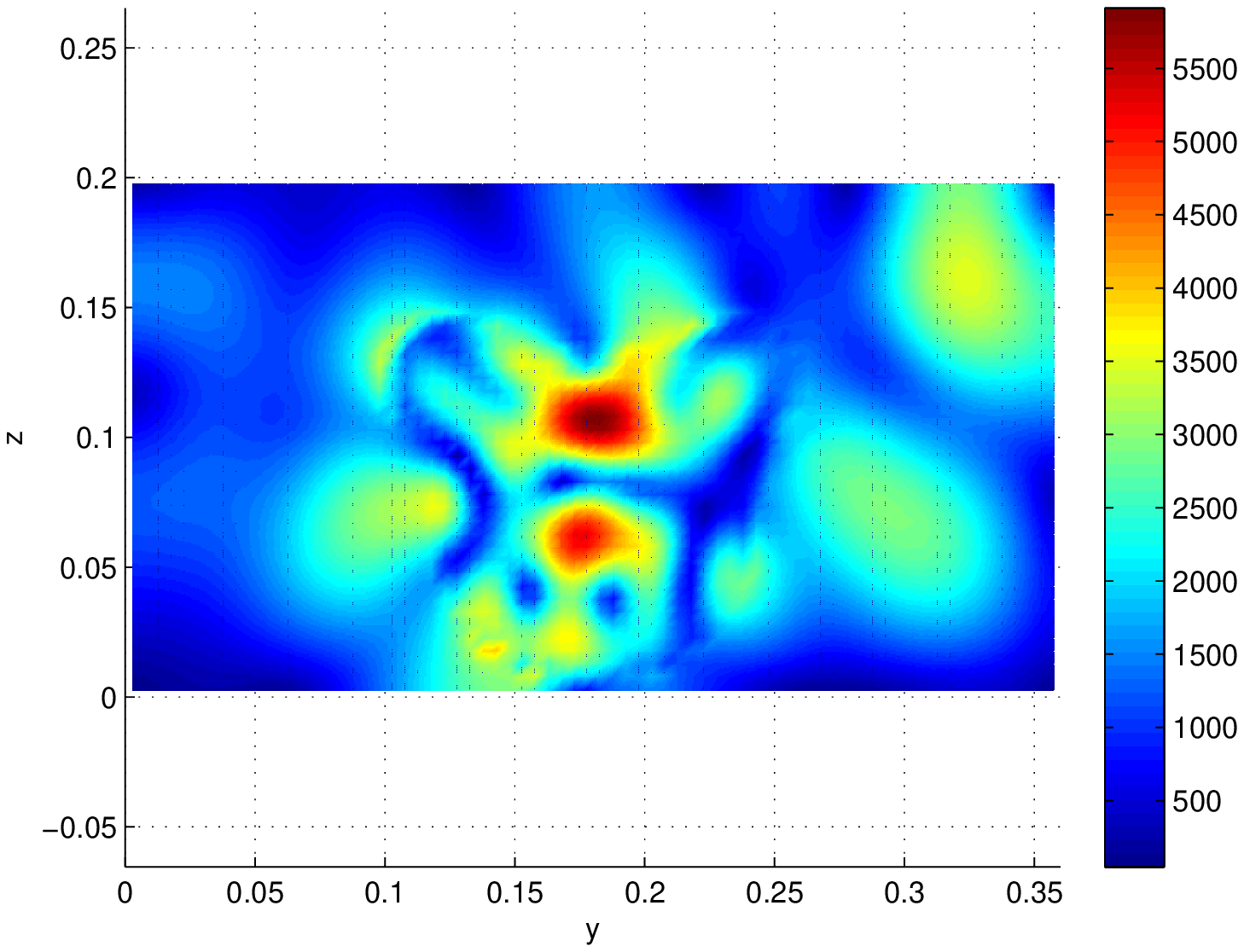} 
  \includegraphics[width=0.49\textwidth]{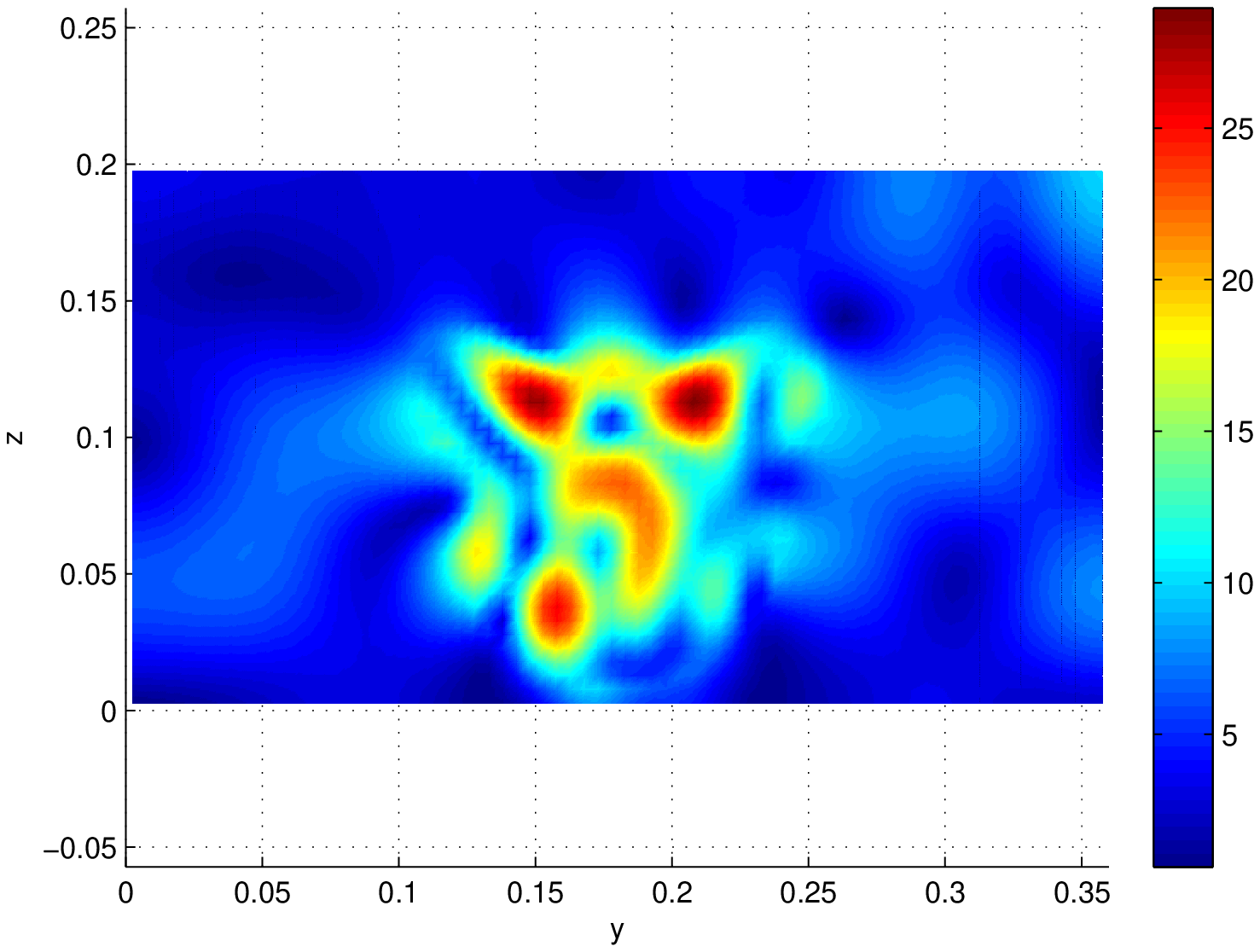}
  \includegraphics[width=0.49\textwidth]{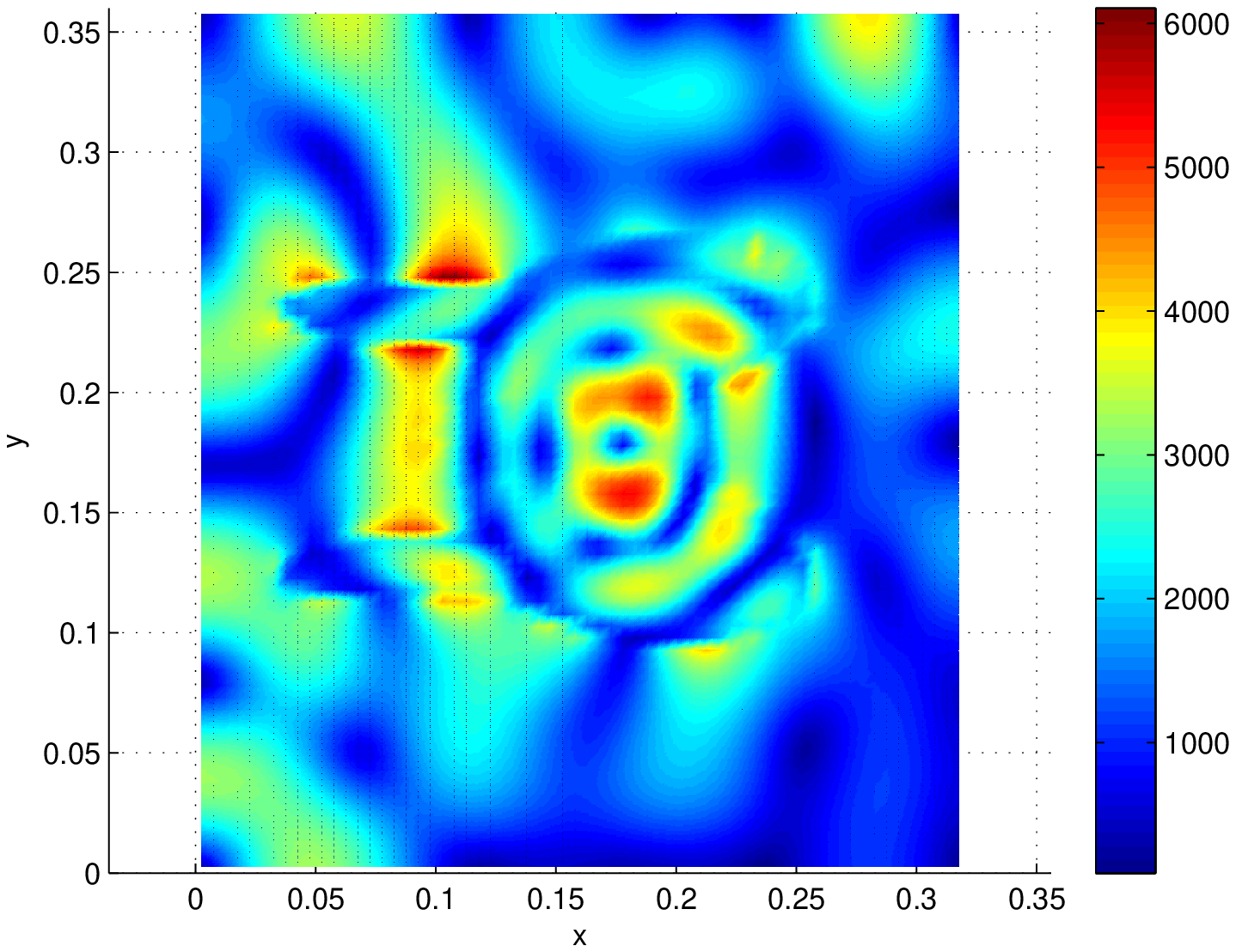} 
  \includegraphics[width=0.49\textwidth]{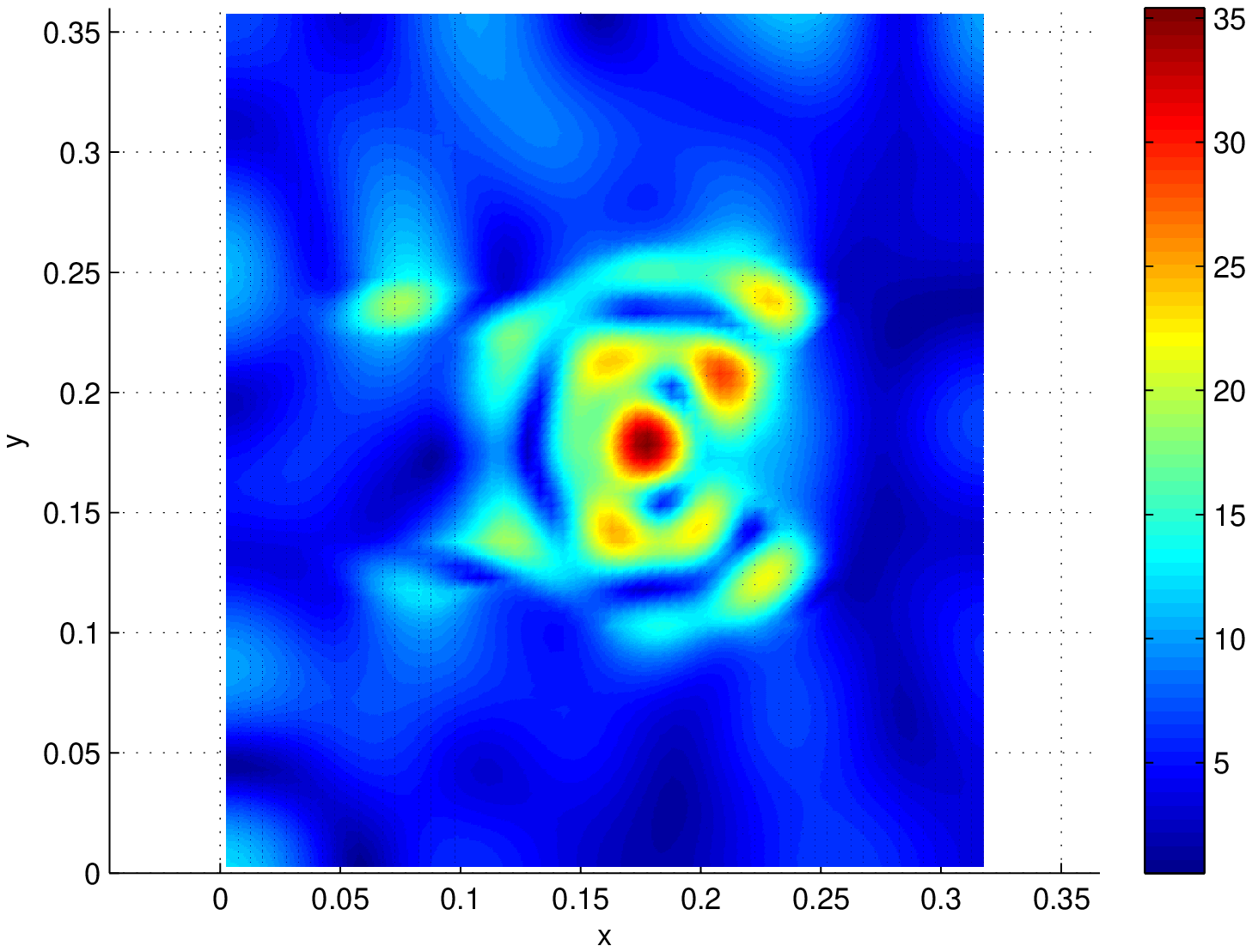}
  \caption{Cross sections of the electric and magnetic field intensity.}
  \label{sol2d}
\end{figure}

We decompose the microwave oven into $2\times 2\times 2=8$ subdomains
of equal size on a grid with mesh size $h=0.005$, which allows us to
solve this problem on a PC, where a direct factorization would not
have been possible any more. The real part of the magnetic field of
the solution is shown in Figure \ref{Chicken} on the right, and the
intensity (Euclidian norm) of the electric and magnetic field in the
oven are shown in Figure \ref{sol3d} in three dimensions.
Two-dimensional cross sections of the solution are shown in Figure
\ref{sol2d},
where we show in each row on the left the electric and on the right
the magnetic field intensity. One can see from these computational
experiments why a turntable is so important in a microwave oven: there
are hot spots, where the intensity of the standing wave is high in the
chicken, and other areas, where there is very little heating happening.
Only a turntable can lead to an approximately even heating of the
chicken.

\section{Conclusions}\label{Conclusions}  
We have shown that for Maxwell's equations, a classical Schwarz
algorithm using characteristic Dirichlet transmission conditions
between subdomains has the same convergence behavior as a simple
optimized Schwarz method applied to the Helmholtz equation, with a low
frequency approximation of the optimal transmission conditions. This
relation allowed us to develop easily an entire hierarchy of
optimized overlapping and non-overlapping Schwarz methods with better
transmission conditions than the characteristic ones for Maxwell's
equations. We illustrated with numerical experiments that the new
algorithms converge much more rapidly than the classical one, and that
such algorithms can be effectively used to compute an approximate
solution for a large scale application. This latter problem contains a
positive conductivity, variable coefficients and multiple subdomains,
a case which is not covered by our current analysis. Nevertheless, the
algorithm performs well with the coefficients derived from the zero
conductivity, constant coefficient case. We are currently studying the
optimization problem with non-zero conductivity, for which the
equivalence with the Helmholtz equation does not hold any longer.

The equivalence between systems and scalar equations has already been
instrumental for the development of optimized Schwarz algorithms for
the Cauchy-Riemann equations, and will almost certainly play an
important role for other cases. For example, it was observed in
\cite{Dolean:2004:CAS} that for Euler's equation, the classical
Schwarz algorithm with characteristic information exchange at the
interfaces is convergent, even without overlap. To relate systems of
partial differential equations to scalar ones, the algebraic tool of
the Smith factorization \cite{Wloka:1995:BVP} has proved to be useful,
see \cite{Dolean:2008:DAN}.  

\bibliographystyle{abbrv}  
\bibliography{ddm,paper} 

\begin{thebibliography}{10}

\bibitem{Alonso:2006:NSM}
A.~Alonso-Rodriguez and L.~Gerardo-Giorda.
\newblock New non-overlapping domain decomposition methods for the
  time-harmonic {M}axwell system.
\newblock {\em SIAM J. Sci. Comp.}, 28(1):102--122, 2006.

\bibitem{Serre:2007:MDH}
S.~Benzoni-Gavage and D.~Serre.
\newblock {\em Multi-dimensional hyperbolic partial differential equations:
  First-order Systems and Applications}.
\newblock Oxford Mathematical Monographs. Oxford Science Publications, 2007.

\bibitem{Chan:1994:DOA}
T.~F. Chan and T.~P. Mathew.
\newblock Domain decomposition algorithms.
\newblock In {\em Acta Numerica 1994}, pages 61--143. Cambridge University
  Press, 1994.

\bibitem{Charton:1991:MDD}
P.~Charton, F.~Nataf, and F.~Rogier.
\newblock M{\'e}thode de d{\'e}composition de domaine pour l'\'equation
  d'advection-diffusion.
\newblock {\em C. R. Acad. Sci.}, 313(9):623--626, 1991.

\bibitem{Chevalier:1998:MNT}
P.~Chevalier.
\newblock {\em M\'ethodes num\'eriques pour les tubes hyperfr\'equences.
  R\'esolution par d\'ecomposition de domaine}.
\newblock PhD thesis, Universit\'e Paris VI, 1998.

\bibitem{Chevalier:1998:SMO}
P.~Chevalier and F.~Nataf.
\newblock Symmetrized method with optimized second-order conditions for the
  {H}elmholtz equation.
\newblock In {\em Domain decomposition methods, 10 (Boulder, CO, 1997)}, pages
  400--407. Amer. Math. Soc., Providence, RI, 1998.

\bibitem{Clerc:1998:NOS}
S.~Clerc.
\newblock Non-overlapping {S}chwarz method for systems of first order
  equations.
\newblock {\em Cont. Math.}, 218:408--416, 1998.

\bibitem{Collino:1997:NIC}
P.~Collino, G.~Delbue, P.~Joly, and A.~Piacentini.
\newblock A new interface condition in the non-overlapping domain
  decomposition.
\newblock {\em Comput. Methods Appl. Mech. Engrg.}, 148:195--207, 1997.

\bibitem{Deng:1997:AND}
Q.~Deng.
\newblock An analysis for a nonoverlapping domain decomposition iterative
  procedure.
\newblock {\em SIAM J. Sci. Comput.}, 18:1517--1525, 1997.

\bibitem{Despres:1990:DDP}
B.~Despr{\'e}s.
\newblock D\'ecomposition de domaine et probl\`eme de {H}elmholtz.
\newblock {\em C.R. Acad. Sci. Paris}, 1(6):313--316, 1990.

\bibitem{Despres:1993:DDM}
B.~Despr{\'e}s.
\newblock Domain decomposition method and the {H}elmholtz problem.{I}{I}.
\newblock In {\em Second International Conference on Mathematical and Numerical
  Aspects of Wave Propagation (Newark, DE, 1993)}, pages 197--206,
  Philadelphia, PA, 1993. {SIAM}.

\bibitem{Despres:1992:ADD}
B.~Despr{\'e}s, P.~Joly, and J.~E. Roberts.
\newblock A domain decomposition method for the harmonic {M}axwell equations.
\newblock In {\em Iterative methods in linear algebra (Brussels, 1991)}, pages
  475--484, Amsterdam, 1992. North-Holland.

\bibitem{Dolean:2007:WCS}
V.~Dolean and M.~J. Gander.
\newblock Why classical {S}chwarz methods applied to hyperbolic systems
  converge even without overlap.
\newblock In {\em Seventeenth International Conference on Domain Decomposition
  Methods}, 2007.

\bibitem{Dolean:2005:IFV}
V.~Dolean and S.~Lanteri.
\newblock An implicite finite volume time-domain method on unstructured meshes
  for {M}axwell equations in three dimensions.
\newblock Technical Report 5767, INRIA, 2005.

\bibitem{Dolean:2002:CIC}
V.~Dolean, S.~Lanteri, and F.~Nataf.
\newblock Construction of interface conditions for solving compressible {E}uler
  equations by non-overlapping domain decomposition methods.
\newblock {\em Int. J. Numer. Meth. Fluids}, 40:1485--1492, 2002.

\bibitem{Dolean:2004:CAS}
V.~Dolean, S.~Lanteri, and F.~Nataf.
\newblock Convergence analysis of a {S}chwarz type domain decomposition method
  for the solution of the {E}uler equations.
\newblock {\em Appl. Num. Math.}, 49:153--186, 2004.

\bibitem{Dolean:2008:DAN}
V.~Dolean, F.~Nataf, and G.~Rapin.
\newblock Deriving a new domain decomposition method for the {S}tokes equations
  using the {S}mith factorization.
\newblock Technical report, http://hal.archives-ouvertes.fr/hal-00110609/fr/,
  2007.

\bibitem{Enquist:1998:ABC}
B.~Engquist and H.-K. Zhao.
\newblock Absorbing boundary conditions for domain decomposition.
\newblock {\em Appl. Numer. Math.}, 27(4):341--365, 1998.

\bibitem{Faccioli:1996:SDM}
E.~Faccioli, F.~Maggio, a.~Quarteroni, and A.~Tagliani.
\newblock Spectral domain decomposition methods for the solution of acoustic
  and elastic wave propagation.
\newblock {\em Geophysics}, 61:1160--1174, 1996.

\bibitem{Faccioli:1997:EWP}
E.~Faccioli, F.~Maggio, A.~Quarteroni, and A.~Tagliani.
\newblock 2d and 3d elastic wave propagation by pseudo-spectral domain
  decomposition method.
\newblock {\em Journal of Seismology}, 1:237--251, 1997.

\bibitem{Gander:2006:OSM}
M.~J. Gander.
\newblock Optimized {S}chwarz methods.
\newblock {\em SIAM J. Numer. Anal.}, 44(2):699--731, 2006.

\bibitem{Gander:2003:MRO}
M.~J. Gander and L.~Halpern.
\newblock M\'ethodes de relaxation d'ondes pour l'\'equation de la chaleur en
  dimension 1.
\newblock {\em C.R. Acad. Sci. Paris, S\'erie I}, 336(6):519--524, 2003.

\bibitem{gander-etal:06}
M.~J. Gander, L.~Halpern, and F.~Magoul\`es.
\newblock An {O}ptimized {S}chwarz {M}ethod with two-sided {R}obin transmission
  conditions for the {H}elmholtz {E}quation.
\newblock {\em Int. J. Numer. Meth. Fluids}, 2007.
\newblock in press.

\bibitem{Gander:2003:OSWW}
M.~J. Gander, L.~Halpern, and F.~Nataf.
\newblock Optimal {S}chwarz waveform relaxation for the one dimensional wave
  equation.
\newblock {\em SIAM Journal of Numerical Analysis}, 41(5):1643--1681, 2003.

\bibitem{Gander:2001:OSH}
M.~J. Gander, F.~Magoul\`es, and F.~Nataf.
\newblock Optimized {S}chwarz methods without overlap for the {H}elmholtz
  equation.
\newblock {\em SIAM J. Sci. Comput.}, 24(1):38--60, 2002.

\bibitem{Hagstrom:2007:RBC}
T.~Hagstrom and S.~Lau.
\newblock Radiation boundary conditions for {M}axwell's equations: a review of
  accurate time-domain formulations.
\newblock {\em J. Comput. Math}, 25(3):305--336, 2007.

\bibitem{Hagstrom:1988:NED}
T.~Hagstrom, R.~P. Tewarson, and A.~Jazcilevich.
\newblock Numerical experiments on a domain decomposition algorithm for
  nonlinear elliptic boundary value problems.
\newblock {\em Appl. Math. Lett.}, 1(3), 1988.

\bibitem{Japhet:2000:OO2}
C.~Japhet, F.~Nataf, and F.~Rogier.
\newblock The optimized order 2 method. application to convection-diffusion
  problems.
\newblock {\em Future Generation Computer Systems FUTURE}, 18, 2001.

\bibitem{Lions:1990:SAM}
P.-L. Lions.
\newblock On the {S}chwarz alternating method. {III:} a variant for
  nonoverlapping subdomains.
\newblock In T.~F. Chan, R.~Glowinski, J.~P{\'e}riaux, and O.~Widlund, editors,
  {\em Third International Symposium on Domain Decomposition Methods for
  Partial Differential Equations , held in Houston, Texas, March 20-22, 1989},
  Philadelphia, PA, 1990. SIAM.

\bibitem{Nataf:1995:FCD}
F.~Nataf and F.~Rogier.
\newblock Factorization of the convection-diffusion operator and the {S}chwarz
  algorithm.
\newblock {\em $M^3AS$}, 5(1):67--93, 1995.

\bibitem{Nedelec:2001:AEE}
J.-C. Nedelec.
\newblock {\em Acoustic and electromagnetic equations. Integral representations
  for harmonic problems.}
\newblock Applied Mathematical Sciences, 144. Springer Verlag, 2001.

\bibitem{Quarteroni:1990:DCL}
A.~Quarteroni.
\newblock Domain decomposition methods for systems of conservation laws:
  spectral collocation approximation.
\newblock {\em SIAM J. Sci. Stat. Comput.}, 11:1029--1052, 1990.

\bibitem{Quarteroni:1996:HDD}
A.~Quarteroni and L.~Stolcis.
\newblock Homogeneous and heterogeneous domain decomposition methods for
  compressible flow at high {R}eynolds numbers.
\newblock Technical Report~33, CRS4, 1996.

\bibitem{Quarteroni:1999:DDM}
A.~Quarteroni and A.~Valli.
\newblock {\em Domain Decomposition Methods for Partial Differential
  Equations}.
\newblock Oxford Science Publications, 1999.

\bibitem{Smith:1996:DPM}
B.~F. Smith, P.~E. Bj{\o}rstad, and W.~Gropp.
\newblock {\em Domain Decomposition: Parallel Multilevel Methods for Elliptic
  Partial Differential Equations}.
\newblock Cambridge University Press, 1996.

\bibitem{Sofronov:1998:NIO}
I.~Sofronov.
\newblock Nonreflecting inflow and outflow in a wind tunnel for transonic
  time-accurate simulation.
\newblock {\em J. Math. Anal. Appl.}, 221(1) 92, 1998.

\bibitem{Sun:1996:OSA}
H.~Sun and W.-P. Tang.
\newblock An overdetermined {S}chwarz alternating method.
\newblock {\em SIAM Journal on Scientific Computing}, 17(4):884--905, Jul.
  1996.

\bibitem{Tang:1992:GSS}
W.~P. Tang.
\newblock Generalized {S}chwarz splittings.
\newblock {\em SIAM J. Sci. Stat. Comp.}, 13(2):573--595, 1992.

\bibitem{Toselli:2000:OSM}
A.~Toselli.
\newblock Overlapping {S}chwarz methods for {M}axwell's equations in three
  dimensions.
\newblock {\em Numer. Math.}, 86(4):733--752, 2000.

\bibitem{Toselli:2004:DDM}
A.~Toselli and O.~Widlund.
\newblock {\em Domain Decomposition Methods - Algorithms and Theory}, volume~34
  of {\em Springer Series in Computational Mathematics}.
\newblock Springer, 2004.

\bibitem{Wloka:1995:BVP}
J.~Wloka, B.~Rowley, and B.~Lawruk.
\newblock {\em Boundary Value Problems for Elliptic Systems}.
\newblock Cambridge University Press, 1995.

\bibitem{Xu:1992:IMS}
J.~Xu.
\newblock Iterative methods by space decomposition and subspace correction.
\newblock {\em SIAM Review}, 34(4):581--613, December 1992.

\bibitem{Xu:1998:NDD}
J.~Xu and J.~Zou.
\newblock Some nonoverlapping domain decomposition methods.
\newblock {\em SIAM Review}, 40:857--914, 1998.

\end{thebibliography}

\end{document}